\documentclass[11pt,reqno]{amsart}
\usepackage{graphicx}
\usepackage{amsmath, amssymb}
\usepackage{verbatim}
\usepackage{color}

\newtheorem{conj.}[thm]{Conjecture}

\theoremstyle{definition}

\theoremstyle{remark}

\numberwithin{equation}{section}

\begin{document}
\title[Quaternion Linear Canonical Wavelet Transform and The Corresponding Uncertainty Inequalities]
{Quaternion Linear Canonical Wavelet Transform and The Corresponding Uncertainty Inequalities}

\author[Aajaz A. Teali]{Aajaz A. Teali}

\address{ Department of Mathematics, University of Kashmir, South Campus, Anantnag-192101 Jammu and Kashmir, India.}
\email{aajaz.math@gmail.com}


\begin{abstract}
 The linear canonical wavelet transform has been shown to be a valuable and powerful time-frequency analyzing tool for optics and signal processing. In this article, we propose a novel transform called quaternion linear canonical wavelet transform which is designed to represent two dimensional quaternion-valued  signals at different  scales, locations and orientations. The proposed transform not only inherits the features of quaternion wavelet transform but also has the capability of signal representation in quaternion linear canonical domain.  We investigate the fundamental properties of quaternion linear canonical wavelet transform including Parseval's formula, energy conservation, inversion formula, and characterization of its range  using the machinery of quaternion linear canonical transform and its convolution. We conclude our investigation by deriving  an analogue of the classical Heisenberg-Pauli-Weyl uncertainty inequality and the associated  logarithmic and local versions for the quaternion linear canonical wavelet transform.
\end{abstract}

\parindent=0mm \vspace{.0in}
\maketitle
\section{Introduction}

\parindent=0mm \vspace{.0in}
Wavelet transforms serve as an important and powerful analyzing tool for time-frequency analysis and  have been applied in a number of fields including signal processing, image processing, sampling theory, differential and integral equations, quantum mechanics and medicine. However, the signal analysis capability of the wavelet transform is limited in the time-frequency plane as each wavelet component is actually a differently scaled bandpass filter in the frequency domain and hence, it does not serve as an efficient tool for processing the higher dimensional signal whose energy is not well concentrated in the frequency domain. One of the examples of such signal is chirp like signals \cite{DS}. Therefore,  in order to obtain joint signal representations in both time and frequency domains, the linear canonical wavelet transform (LCWT) has been introduced in the context of time-frequency analysis. The  LCWT inherits the excellent mathematical properties of wavelet transform and linear canonical transform along with some fascinating properties of its own. In recent years, this transform has been paid a considerable amount of attention, resulting in many applications in the areas of optics, quantum mechanics, pattern recognition and signal processing \cite{Wei, Wang, Yong, Prasad, Yong2}.

\parindent=8mm \vspace{.1in}
On the other hand, considerable attention has been paid for the representation of signals in quaternion domains as quaternion algebra  is the closest in its mathematical properties to the familiar system of the real and complex numbers. The quaternion algebra offers a simple and profound representation of signals wherein several components are to be controlled simultaneously. The development of integral transforms for quaternion valued signals has found numerous applications in 3D computer graphics, aerospace engineering,  artificial intelligence and colour image processing. The extension of classical  wavelet transform to quaternion algebra has been introduced in \cite{He} and \cite{Zhao}, they also demonstrated their various properties. In \cite{Tra}, Traversoni proposed a discrete quaternion wavelet transform, the application of which can be found in \cite{Bayro,Maw} and \cite{Fletcher}.

\parindent=8mm \vspace{.1in}
In the recent years, some authors have generalized the linear canonical transform to quaternion-valued signals, known as the quaternionic linear
canonical transform (QLCT). The QLCT was firstly studied in \cite{Kou2} including prolate spheroidal wave signals and uncertainty principles \cite{Kou}. Some useful properties of the QLCT such as linearity, reconstruction formula, continuity, boundedness, positivity inversion formula and the uncertainty principle were established in \cite{Asym,MR,Zhang,Zhou}. An application of the QLCT to study of generalized swept-frequency filters was introduced in  \cite{BahApp}. Because of the non-commutative property of multiplication of quaternions, there are mainly three various types of 2D quaternion linear canonical transform (QLCTs): two-sided QLCTs, left-sided QLCTs and right-sided QLCTs (refer to \cite{Kou}). Based on the (two-sided) QLCT \cite{Kou2}, the quaternion windowed linear canonical transform of 2D quaternionic signals has been introduced by Wen-Biao Gao and Bing-Zhao Li in \cite{Wen}. It can reveal the local QLCT-frequency contents and enjoys high concentrations and eliminates the cross term, but it has limitation of having fixed window localization.

\parindent=8mm \vspace{.1in}
The objective of this paper is to eliminate the redundancy of the quaternion wavelet transform(QWT) and quaternion linear canonical transform(QLCT) by introducing the quaternion linear canonical wavelet transform(QLCWT) i.e. the quaternion version of linear canonical wavelet transform. QLCWT generalizes the definition of the QWT in time-QLCD-frequency plane by using the modified quaternion linear canonical convolution. The proposed QLCWT not only inherits the features of quaternion wavelet transform but also has the capability of signal representation in quaternion domain. Besides, it has explicit physical interpretation and low complexity. It is hoped that this transform might be useful in three dimensional color field processing, space color video processing, crystallography, aerospace engineering, oil exploration and for the solution of many types of quaternionic differential equations.

\parindent=8mm \vspace{.1in}
The article is organized as follows: We begin in Section 2 by presenting the notation, quaternion algebra and wavelet theory needed to understand and place our results in context. In Section 3, we introduce the concept of quaternion linear canonical wavelet transform and obtain the expected properties of the extended linear canonical wavelet transform  including Parseval's formula, energy conservation, inversion formula, and characterization of its range. The well known Heisenberg-Pauli-Weyl inequality, logarithmic and local uncertainty principle are generalized in the quaternion Linear canonical domains in Section 4.

\section{Preliminaries}

\subsection{ Quaternion Algebra}
\parindent=0mm \vspace{.1in}
The theory of quaternions was initiated by the Irish mathematician Sir W.R. Hamilton in 1843 and is denoted by $\mathbb H$ in his honour.
The quaternion algebra provides an extension of the complex number system to an associative non-commutative four-dimensional algebra.
The quaternion algebra $\mathbb H$ over $\mathbb R$ is given by
\begin{align*}
\mathbb H=\Big\{f=a_{0}+i\,a_{1}+j\,a_{2}+k\,a_{3}\,:\,a_{0},a_{1},a_{2},a_{3}\in\mathbb R\Big\},
\end{align*}
where $i,j,k$ denote the three imaginary units, obeying the Hamilton's multiplication rules
\begin{align*}
ij=k=-ji,~jk=i=-kj\,,\,ki=j=-ik\,,\, {\text {and}}\,\, {i}^{2}={j}^{2}={k}^{2}={ijk}=-1.
\end{align*}
For  quaternions ${f}_{1}=a_{0}+i\,a_{1}+j\,a_{2}+k\,a_{3}$ and $ {f}_{2}=b_{0}+i\,b_{1}+ j\,b_{2}+k\,b_{3}$, the addition  is defined componentwise
and the multiplication is defined as
\begin{align*}
{f}_{1} {f}_{2}&=(a_{0}b_{0}-a_{1}b_{1}-a_{2}b_{2}-a_{3}b_{3})+i\,(a_{1}b_{0}+a_{0}b_{1}+a_{2}b_{3}-a_{3}b_{2})\\
&~~~~~~~~~~~~~~+ j(a_{0}b_{2}+a_{2}b_{0}+a_{3}b_{1}-a_{1}b_{3})+k\,(a_{0}b_{3}+a_{3}b_{0}+a_{1}b_{2}-a_{2}b_{1}).
\end{align*}
The conjugate and norm of a quaternion  $f=a_{0}+i\,a_{1}+j\,a_{2}+k\,a_{3},$ are given by $\overline{f}=a_{0}-i\,a_{1}-j\,a_{2}-k\,a_{3}$ and ${\| f\|}_{\mathbb H}=\sqrt{{a_{0}}^{2}+{a_{1}}^{2}+{a_{2}}^{2}+{a_{3}}^{2}}$, respectively. We also note that an  arbitrary quaternion $h$ can be represented by two complex numbers  as $f=(a_{0}+i\,a_{1})+j\,(a_{2}-i\,a_{3})=f_1+j\,f_2$, where $f_1,f_2\in\mathbb C$, and hence, $\overline{f}=\overline {f_1}-j\,f_2$, with $\overline{f_1}$ denoting the complex conjugate of $f_1$. Moreover, the inner product of any two quaternions  $f=f_{1}+\,jf_{2},$ and $g=g_{1}+j\,g_{2}$ in ${\mathbb H}$ is defined by

\begin{align*}
{\big\langle f,\,g\big\rangle}_{\mathbb H}&= f\overline{g}=(f_{1}{\overline g_{1}}+{\overline f_{2}}g_{2})+j(f_{2}{\overline g_{1}}-{\overline f_{1}}g_{2}).
\end{align*}

By virtue of the complex domain representation, a quaternion-valued function $f:\mathbb R^2 \to \mathbb H$ can be decomposed as $f(x)=f_{1}+j\,f_{2}$,  where $f_{1},f_{2}$ are both complex valued functions.

\parindent=8mm \vspace{.1in}
Let us denote $L^2(\mathbb R^2,\mathbb H)$, the space of all quaternion valued functions $f$ satisfying
\begin{align*}
 {\big\|f\big\|}_{2}=\left\{\int_{\mathbb R^{2}}\Big(|f_{1}({\bf x})|^{2}+|f_{2}({\bf x})|^{2}\Big)\,d{\bf x}\right\}^{1/2}< \infty.
\end{align*}
The norm on $L^2(\mathbb R^2,\mathbb H)$ is obtained from the inner product of the quaternion valued functions $f=f_{1}+j\,f_{2},$ and $g=g_{1}+j\,g_{2}$  as
\begin{align*}
{\big\langle f,\,g\big\rangle}_{2}&= \int_{\mathbb R^{2}}{\big\langle f,\,g\big\rangle}_{\mathbb H}\,d{\bf x}\\
&=\int_{\mathbb R^{2}}\bigg\{\Big(f_{1}({\bf x})\,\overline{g_{1}}({\bf x})+\overline{f_{2}}({\bf x})\,g_{2}({\bf x})\Big)+j \Big(f_{2}({\bf x})\,\overline{g_{1}}({\bf x})
-\overline{f_{1}}({\bf x})\,g_{2}({\bf x})\Big)\bigg\}\,d{\bf x}.
\end{align*}
An easy computation shows that $L^2(\mathbb R^2,\mathbb H)$ equipped with above defined inner product is a Hilbert space.

\parindent=0mm \vspace{.1in}
{\bf Definition 2.1.}  For any quaternion valued function $f\in L^1(\mathbb R^2,\mathbb H)\cap L^2(\mathbb R^2,\mathbb H)$, the two-sided quaternion Fourier transform (QFT) s denoted by ${\mathcal F_{q}}$ and is given by
\begin{align*}
{\mathcal F_{q}}\big[f({\bf x})\big]({\bf w}) =\hat f({\bf w})= \int_{\mathbb R^2}e^{- i x_{1}  w_{1}} f({\bf x})\,e^{- j x_{2}  w_{2}}\,d{\bf x},\tag{2.1}
\end{align*}
where ${\bf x}=(x_{1},x_{2}),\, {\bf w}=(w_{1},w_{2})$ and the quaternion exponential $e^{- i x_{1}  w_{1}}$ and $e^{- j x_{2}  w_{2}}$ are the quaternion Fourier kernels. The corresponding inversion formula is given by
\begin{align*}
f({\bf x}) &=\dfrac{1}{(2\pi)^2}\int_{\mathbb R^2}e^{i x_{1} w_{1}} f({\bf x})\,e^{j x_{2} w_{2}}\,d{\bf w}.\tag{2.2}
\end{align*}

\subsection{ Quaternion Linear Canonical Transform }

\parindent=0mm \vspace{.1in}
Due to non-commutativity of quaternion multiplication, there are three types of the quaternion offset linear canonical transform, the left-sided QLCT, the right-sided QLCT, and two-sided QLCT.  In this paper, we will extend the theory of two-sided QLCT introduced by  \cite{Kou2} as follows:

\parindent=0mm \vspace{.1in}
{\bf Definition 2.2. } Let $A_s =\left[\begin{array}{cc}a_s & b_s \\c_s & d_s \\ \end{array}\right] $, be a matrix parameter such that $a_s, b_s, c_s, d_s \in \mathbb R$ and $ a_sd_s-b_sc_s=1,$ for $s=1,2.$ The two-sided quaternion linear canonical transform of any quaternion valued function $f\in L^2(\mathbb R^2,\mathbb H)$, is given by

\begin{align*}
\mathcal L^{\mathbb H}_{A_1,A_2}\big[f({\bf x})\big]({\bf w})=\left\{\begin{array}{cc}\int_{\mathbb R^2} K_{A_1}^i(x_1,w_1) \,f({\bf x})\,K_{A_2}^j(x_2,w_2)d{\bf x},& b_1b_2\neq0\\ \sqrt{d_1d_2} e^{\frac{ic_1d_1}{2}w_1^2} f(d_1w_1,d_2w_2) e^{\frac{jc_2d_2}{2}w_2^2}, &b_1b_2=0\\
\end{array}\right.
\end{align*}
where ${\bf x}=(x_{1},x_{2}),\, {\bf w}=(w_{1},w_{2})$ and for $b_1b_2\neq0$ ,the quaternion kernels $K_{A_1}^i(x_1,w_1)$ and $K_{A_2}^j(x_2,w_2)$ are respectively given by
\begin{align*}
K_{A_1}^i(x_1,w_1)=\dfrac{1}{\sqrt{2\pi b_1}} \exp\left\{\frac{i}{2b_1}\big[a_1x_1^2-2x_1w_1+d_1w_1^2-\frac{\pi b_1}{2}\big]\right\}
\end{align*}
and

\begin{align*}
K_{A_2}^j(x_2,w_2)=\dfrac{1}{\sqrt{2\pi b_2}} \exp\left\{\frac{j}{2b_2}\big[a_2x_2^2-2x_2w_2+d_2w_2^2-\frac{\pi b_2}{2}\big]\right\},
\end{align*}

\parindent=0mm \vspace{.1in}
The corresponding inversion formula for two-sided QOLCT is given by
\begin{align*}
f({\bf x}) &=\int_{\mathbb R^2}\overline{K_{A_1}^i(x_1,w_1)} \, \mathcal L^{\mathbb H}_{A_1,A_2}\big[f\big]({\bf w})\,\overline{K_{A_2}^j(x_2,w_2)} \,d{\bf w}.\tag{2.3}
\end{align*}

\parindent=0mm \vspace{.0in}
Two quaternion functions $f,\,g \in L^2(\mathbb R^2,\mathbb H)$ are related to their two-sided QOLCT via the Parseval formula, given as

\begin{align*}
\Big\langle \mathcal L^{\mathbb H}_{A_1,A_2} [f],\,\mathcal L^{\mathbb H}_{A_1,A_2} [g]\Big\rangle_{L^2(\mathbb R^2,\mathbb H)}=\big\langle f,g \big\rangle_{\mathbb H}.\,\tag{2.4}
\end{align*}
In particular, $\Big\|\mathcal L^{\mathbb H}_{A_1,A_2}\big[f({\bf x})\big]({\bf w})\Big\|_{L^2(\mathbb R^2,\mathbb H)}=\left\|f\right\|_{L^2(\mathbb R^2,\mathbb H)}.$

{\bf Definition 2.3.} The set
$$ \mathcal G = \mathbb R^{+} \times \mathbb R^{2}\times \text{SO}(2)=\left\{(a,{\bf y}, r_\theta);\,a\in \mathbb R^{+},\, {\bf y}\in\mathbb R^{2},\, r_\theta \in\text{SO}(2) \right\}$$

endowed with the operation
\begin{align*}
\big (a,{\bf y}, r_\theta \big )\odot \big (a^{\prime},{\bf y}^{\prime},r_{\theta}^{\prime}\big )= \big (aa^{\prime},{\bf y}+ar_{\theta} {\bf y}^{\prime},r_{\theta+\theta^{\prime}}\big),
\end{align*}
forms a group called a {\it similitude group} on $\mathbb R^{2}$ associated with wavelets, where SO(2) is the special orthogonal group of rotation in $\mathbb R^{2}$.  The left  Haar measures on $\mathcal G$ is given by  $d\eta={da\,d{\bf y}\, d\theta}/a^{3}$ . For more see \cite{Antoine}

\section{Quaternion Linear Canonical Wavelet Transform}

\parindent=0mm \vspace{.1in}
In this section, we shall first introduce a new convolution structure for two dimensional linear canonical transform and obtain the corresponding convolution theorem. Then, we shall characterize the admissibility condition in terms of two-sided QLCT and define the quaternion linear canonical wavelet transform in terms of admissible canonical wavelets based on convolution operator in LCT domain and investigate its fundamental properties.

\parindent=0mm \vspace{.1in}
Analogues to Wei et al \cite{MultiC}, the generalized convolution theorem for two dimensional linear canonical transform, based on generalized translation is derived as follows.

\parindent=0mm \vspace{.1in}
The ${\bf y}$-generalized translation of signal $\psi({\bf x})$ is denoted by $\psi({\bf x}\Theta{\bf y})$ is given by \cite{Marks}

\begin{align*}
\psi({\bf x}\Theta{\bf y})= \int \rho({\bf w}) \Psi({\bf w})K({\bf x},{\bf w}) K^{*}({\bf y},{\bf w})d{\bf w},\,\tag{3.1}
\end{align*}

where $\Theta$ is argument of function, $\psi({\bf x}\Theta{\bf y})$ is generalized delay operator for the generalized translation, $\rho$ is weight function, $\Psi({\bf w})$ is transformed function $\psi({\bf x})$, and $K({\bf x},{\bf w})$ is kernel of transformation. The corresponding generalized translation in two sided QLC-domain is given by

\begin{align*}
\psi_{A_s}({\bf x}\Theta{\bf y})&= \int_{\mathbb R^2} K_{A_1}^{-i}(y_1,w_1) K_{A_1}^i(x_1,w_1) \mathcal L^{\mathbb H}_{A_1,A_2}\big[\psi({\bf x})\big]({\bf w})K_{A_2}^j(x_2,w_2)K_{A_1}^{-j}(y_2,w_2)d{\bf w}\\
&=\dfrac{1}{(2\pi)^2b_1b_2}\int_{\mathbb R^2}  \exp\left\{\frac{i}{2b_1}\big[a_1(x_1^2-y_1^2)-2(x_1-y_1)w_1\big]\right\} \mathcal L^{\mathbb H}_{A_1,A_2}\big[\psi({\bf x})\big]({\bf w})\\
&\qquad\qquad\qquad\qquad\qquad\qquad\qquad\times\exp\left\{\frac{j}{2b_2}\big[a_2(x_2^2-y_2^2)-2(x_2-y_2)w_2\big]\right\}d{\bf w}\\
&=e^{\frac{ia_1}{2b_1}(x_1^2-y_1^2)}\left\{\dfrac{1}{(2\pi)^2b_1b_2}\int_{\mathbb R^2}  e^{\frac{i(y_1-x_1)w_1}{b_1}} \mathcal L^{\mathbb H}_{A_1,A_2}\big[\psi({\bf x})\big]({\bf w})e^{\frac{j(y_2-x_2)w_2}{b_2}}d{\bf w}\right\}e^{\frac{ja_2}{2b_2}(x_2^2-y_2^2)}.\,\tag{3.2}\\
%
\end{align*}

\parindent=0mm \vspace{.1in}
We now introduce a new convolution structure for two dimensional linear canonical transform and obtain the corresponding convolution theorem.

\parindent=0mm \vspace{.1in}
{\bf Definition 3.1.} For any two signals $f,\,\psi\in L^2(\mathbb R^2,\,\mathbb H),$ we define the generalized convolution operator $\otimes^{A_s}$ by

\begin{align*}
f\otimes^{A_s} \psi({\bf x})= \int_{\mathbb R^2} f({\bf y})\psi_{A_s}({\bf x}\Theta{\bf y}) d{\bf y},\,\tag{3.3}
\end{align*}
where $\psi_{A_s}({\bf x}\Theta{\bf y})$ is given by (3.2), and the corresponding convolution theorem is

\parindent=0mm \vspace{.1in}
{\bf Theorem 3.2.} If $f\otimes^{A_s} \psi({\bf x})$ is defined as (3.3) and $\mathcal L_{A_1,A_2}^{\mathbb H}\big[\psi({\bf x})\big]\in L^2(\mathbb R^2,\,\mathbb R)$ denotes the two dimensional linear canonical transform of a quaternion valued signal $\psi.$ Then we have

\begin{align*}
\mathcal L_{A_1,A_2}^{\mathbb H}\big[f\otimes^{A_s} \psi({\bf x})\big]({\bf w})=\mathcal L_{A_1,A_2}^{\mathbb H}\big[f({\bf x})\big]({\bf w})\cdot\mathcal L_{A_1,A_2}^{\mathbb H}\big[\psi({\bf x})\big]({\bf w}).\,\tag{3.4}
\end{align*}

{\it Proof.} By definition of QLCT, and well known Fubini Theorem, we have

\begin{align*}
&\mathcal L_{A_1,A_2}^{\mathbb H}\big[f\otimes^{A_s} \psi({\bf x})\big]({\bf w})\\
&\qquad= \int_{\mathbb R^2} K_{A_1}^i(x_1,w_1) \,f\otimes^{A_s} \psi({\bf x})\,K_{A_2}^j(x_2,w_2)d{\bf x}\\
&\qquad= \int_{\mathbb R^2} K_{A_1}^i(x_1,w_1) \,\int_{\mathbb R^2} f({\bf y})\psi({\bf x}\Theta{\bf y}) d{\bf y}\,K_{A_2}^j(x_2,w_2)d{\bf x}\\
&\qquad= \int_{\mathbb R^2}\int_{\mathbb R^2} K_{A_1}^i(x_1,w_1) \, f({\bf y})\,\,e^{\frac{ia_1}{2b_1}(x_1^2-y_1^2)}\\
&\qquad\times\left\{\dfrac{1}{(2\pi)^2b_1b_2}\int_{\mathbb R^2}  e^{\frac{i(y_1-x_1)w_1}{b_1}} \mathcal L^{\mathbb H}_{A_1,A_2}\big[\psi({\bf x})\big]({\bf w})e^{\frac{j(y_2-x_2)w_2}{b_2}}d{\bf w}\right\}e^{\frac{ja_2}{2b_2}(x_2^2-y_2^2)}\,\, \,K_{A_2}^j(x_2,w_2)d{\bf x}d{\bf y}\\
&\qquad= \int_{\mathbb R^2}\int_{\mathbb R^2} K_{A_1}^i(x_1,w_1) \, f({\bf y})\,\,e^{\frac{ia_1}{2b_1}(x_1^2-y_1^2)}\\
&\qquad\times\left\{\dfrac{1}{(2\pi)^2b_1b_2}\int_{\mathbb R}  e^{\frac{i(y_1-x_1)w_1}{b_1}}dw_1\,\int_{\mathbb R}e^{\frac{j(y_2-x_2)w_2}{b_2}}dw_2\right\}e^{\frac{ja_2}{2b_2}(x_2^2-y_2^2)}\,\,  \mathcal L^{\mathbb H}_{A_1,A_2}\big[\psi({\bf x})\big]({\bf w})\,K_{A_2}^j(x_2,w_2)\,d{\bf x}d{\bf y}\\
&\qquad= \int_{\mathbb R^2}\int_{\mathbb R^2} K_{A_1}^i(x_1,w_1) \, f({\bf y})\,\,e^{\frac{ia_1}{2b_1}(x_1^2-y_1^2)}\\
&\qquad\times\left\{\dfrac{1}{b_1b_2} b_1\delta(y_1-x_1)\cdot b_2\delta(x_2-y_2) \right\}e^{\frac{ja_2}{2b_2}(x_2^2-y_2^2)}\,\,  \mathcal L^{\mathbb H}_{A_1,A_2}\big[\psi({\bf x})\big]({\bf w})\,K_{A_2}^j(x_2,w_2)\,d{\bf x}d{\bf y}\\
&\qquad= \int_{\mathbb R^2}\int_{\mathbb R^2} K_{A_1}^i(x_1,w_1) \, f({\bf y})\,\,e^{\frac{ia_1}{2b_1}(x_1^2-y_1^2)}\delta(y_1-x_1)\cdot \delta(x_2-y_2)e^{\frac{ja_2}{2b_2}(x_2^2-y_2^2)}\\
&\qquad\qquad\qquad\qquad\qquad\qquad\qquad\qquad\times  \mathcal L^{\mathbb H}_{A_1,A_2}\big[\psi({\bf x})\big]({\bf w})\,K_{A_2}^j(x_2,w_2)\,d{\bf x}d{\bf y}\\
&\qquad= \int_{\mathbb R^2} K_{A_1}^i(x_1,w_1) \, f({\bf y})\, \mathcal L^{\mathbb H}_{A_1,A_2}\big[\psi({\bf x})\big]({\bf w})\,K_{A_2}^j(x_2,w_2)\,d{\bf x}\\
&\qquad= \int_{\mathbb R^2} K_{A_1}^i(x_1,w_1) \, f({\bf y})K_{A_2}^j(x_2,w_2)\,d{\bf x}\cdot \mathcal L^{\mathbb H}_{A_1,A_2}\big[\psi({\bf x})\big]({\bf w})\,\\
&\qquad=\mathcal L_{A_1,A_2}\big[f({\bf x})\big]({\bf w})\cdot\mathcal L_{A_1,A_2}\big[\psi({\bf x})\big]({\bf w}).
\end{align*}
This completes the proof. \quad\fbox

\parindent=0mm \vspace{.1in}

Keeping in view the generalized translation in equation (3.2), we have the following definition of quaternion linear canonical wavelet family of $\psi \in L^2(\mathbb R^2, \,\mathbb H)$.

\parindent=0mm \vspace{.1in}

{\bf Definition 3.3.} For a quaternion signal $\psi\in L^2(\mathbb R^2, \,\mathbb H),$ we define quaternion linear canonical wavelets as
\begin{align*}
U_{a,{\bf y},\theta}\Psi({\bf x})=\Psi_{a,{\bf y},\theta}^{\mathbb H}({\bf x})= a^{-1} e^{\frac{-ia_1}{2b_1}(x_1^2-y_1^2)} \psi\big(r_{-\theta}a^{-1}({\bf x}-{\bf y})\big)\,e^{\frac{-ja_2}{2b_2}(x_2^2-y_2^2)},\,\tag{3.5}
\end{align*}
where $U_{a,{\bf y},\theta}: L^2(\mathbb R^2, \,\mathbb H)\rightarrow L^2(\mathcal G, \,\mathbb H) $ is a unitary operator, $a\in\mathbb R^{+},{\bf y}\in\mathbb R^2,$ and $r_{-\theta}\in $ SO(2), with $r_{-\theta}{\bf x}=(x_1\cos\theta+x_2\sin\theta,\,-x_1\sin\theta+x_2\cos\theta),\,0\le\theta\le2\pi$. The order of the terms in (3.5) is fixed because of the non-commutativity of quaternions.

\parindent=8mm \vspace{.0in}
We now  prove a lemma which offers us the two-sided quaternion LCT of $\Psi_{a,{\bf y},\theta}^{\mathbb H}({\bf x})$, defined in (3.5), in terms of quaternion Fourier transform, which will be fruitful for investigating certain fundamental properties.

\parindent=0mm \vspace{.1in}
{\bf Lemma 3.4.} The two-sided quaternion LCT of $\Psi_{a,{\bf y},\theta}^{\mathbb H}({\bf x})$ , defined in (3.5) is given by

\begin{align*}
\mathcal L^{\mathbb H}_{A_1,A_2}\Big[\Psi_{a,{\bf y},\theta}^{\mathbb H}({\bf x})\Big]({\bf w})&= \dfrac{a}{2\pi \sqrt{b_1b_2}} \exp\left\{\frac{i}{2b_1}\big[d_1w_1^2+y_1^2a_1-2y_1w_1-\frac{\pi b_1}{2}\big]\right\} \mathcal F_q[\psi]\left(r_{-\theta}\,\dfrac{a{\bf w}}{{\bf b}}\right)\\
&\qquad\qquad\qquad\times\exp\left\{\frac{j}{2b_2}\big[d_2w_2^2+y_2^2a_2-2y_2w_2-\frac{\pi b_2}{2}\big]\right\}.\,\tag{3.6}
\end{align*}

\parindent=0mm \vspace{.1in}
{\it Proof.} We have

\begin{align*}
&\mathcal L^{\mathbb H}_{A_1,A_2}\Big[\Psi_{a,{\bf y},\theta}^{\mathbb H}({\bf x})\Big]({\bf w})\\
&\qquad= \dfrac{1}{2\pi a \sqrt{b_1b_2}} \int_{\mathbb R^2}\exp\left\{\frac{i}{2b_1}\big[a_1x_1^2-2x_1w_1+d_1w_1^2-\frac{\pi b_1}{2}\big]\right\}\\
&\qquad\qquad\times e^{\frac{-ia_1}{2b_1}(x_1^2-y_1^2)} \psi\big(r_{-\theta}a^{-1}({\bf x}-{\bf y})\big)\,e^{\frac{-ja_2}{2b_2}(x_2^2-y_2^2)}\exp\left\{\frac{j}{2b_2}\big[a_2x_2^2-2x_2w_2+d_2w_2^2-\frac{\pi b_2}{2}\big]\right\}d{\bf x}\\
&\qquad= \dfrac{1}{2\pi a \sqrt{b_1b_2}} \exp\left\{\frac{i}{2b_1}\big[d_1w_1^2-\frac{\pi b_1+a_1y_1^2}{2}\big]\right\} \int_{\mathbb R^2}e^{\frac{-ix_1w_1}{b_1}} \psi\big(r_{-\theta}a^{-1}({\bf x}-{\bf y})\big)\,e^{\frac{-jx_2w_2}{b_2}}\, d{\bf x}\\
&\qquad\qquad\qquad\qquad\qquad\qquad\qquad\qquad\times\exp\left\{\frac{j}{2b_2}\big[a_2y_2^2+d_2w_2^2-\frac{\pi b_2}{2}\big]\right\}\\
&\qquad= \dfrac{1}{2\pi a \sqrt{b_1b_2}} \exp\left\{\frac{i}{2b_1}\big[d_1w_1^2-\frac{\pi b_1+a_1y_1^2}{2}\big]\right\} \int_{\mathbb R^2}e^{\frac{-i(az_1+y_1)w_1}{b_1}} \psi\big(r_{-\theta}z\big)\,e^{\frac{-j(az_2+y_2)w_2}{b_2}}\, a^2dz\\
&\qquad\qquad\qquad\qquad\qquad\qquad\qquad\qquad\times\exp\left\{\frac{j}{2b_2}\big[a_2y_2^2+d_2w_2^2-\frac{\pi b_2}{2}\big]\right\}\\
&\qquad= \dfrac{a}{2\pi \sqrt{b_1b_2}} \exp\left\{\frac{i}{2b_1}\big[d_1w_1^2-\frac{\pi b_1+a_1y_1^2-2y_1w_1}{2}\big]\right\} \int_{\mathbb R^2}e^{\frac{-i(az_1w_1}{b_1}} \psi\big(r_{-\theta}z\big)\,e^{\frac{-jaz_2w_2}{b_2}}\, dz\\
&\qquad\qquad\qquad\qquad\qquad\qquad\qquad\qquad\times\exp\left\{\frac{j}{2b_2}\big[-2y_2w_2+a_2y_2^2+d_2w_2^2-\frac{\pi b_2}{2}\big]\right\}\\
&\qquad= \dfrac{a}{2\pi \sqrt{b_1b_2}} \exp\left\{\frac{i}{2b_1}\big[d_1w_1^2-\frac{\pi b_1+a_1y_1^2-2y_1w_1}{2}\big]\right\} \,\mathcal F_q[\psi]\left(r_{-\theta}\,\dfrac{a{\bf w}}{{\bf b}}\right)\\
&\qquad\qquad\qquad\qquad\qquad\qquad\qquad\qquad\times\exp\left\{\frac{j}{2b_2}\big[-2y_2w_2+a_2y_2^2+d_2w_2^2-\frac{\pi b_2}{2}\big]\right\}.\\
\end{align*}

\parindent=0mm \vspace{.1in}
{\bf Remark:} From above lemma, we note the following relation
\begin{align*}
&\frac{2\pi \sqrt{b_1b_2}}{a}\mathcal L^{\mathbb H}_{A_1,A_2}\Big[\Psi_{a,{\bf y},\theta}^{\mathbb H}({\bf x})\Big]({\bf b}{\bf w})\\
&\qquad=\exp\left\{\frac{i}{2b_1}\big[d_1(b_1w_1)^2-\frac{\pi b_1+a_1y_1^2-2y_1b_1w_1}{2}\big]\right\}\mathcal F_q[\psi]\big(r_{-\theta}a{\bf w}\big)\\
&\qquad\qquad\qquad\qquad\qquad\qquad\qquad\qquad\times\exp\left\{\frac{j}{2b_2}\big[-2y_2b_2w_2+a_2y_2^2+d_2(b_2w_2)^2-\frac{\pi b_2}{2}\big]\right\}\,\tag{3.7}\\
\end{align*}

In a sequel, we have the following admissibility condition in terms of two-sided quaternion linear canonical transform.

\parindent=0mm \vspace{.0in}
{\bf Definition 3.5(Admissible Canonical Quaternion Wavelet).} A quaternion-valued wavelet $\psi \in L^2\left(\mathbb R^{2},\,\mathbb H\right)$ is admissible, if\\
\begin{align*}
C_\psi=\int_{\mathbb R^{+}}\int_{\text{SO}(2)} \left|\mathcal F_q[\psi]\big(r_{-\theta}a{\bf w}\big)\right|^2\dfrac{dad\theta}{a}= \int_{\mathbb R^{+}\times\text{SO(2)}}\left|\mathcal L^{\mathbb H}_{A_1,A_2}\Big[\Psi_{a,{\bf y},\theta}^{\mathbb H}({\bf x})\Big]({\bf w})\right|^2\,\dfrac{da  d\theta}{a^3}\,\tag{3.8}
\end{align*}
is invertible real constant for a.e. ${\bf w}\in \mathbb R^2$.

\parindent=8mm \vspace{.1in}
We are now ready to define the quaternion linear canonical wavelet transform of two dimensional quaternion valued signals.

\parindent=0mm \vspace{.0in}
{\bf Definition 3.6.} The quaternion linear canonical wavelet transform based on 2D-LCT convolution, of a quaternion-valued function $f \in L^2\left(\mathbb R^{2},\,\mathbb H\right)$ with respect to an admissible wavelet $\psi \in L^2\left(\mathbb R^{2},\,\mathbb H\right)$ is defined by

\begin{align*}\label{3.9}
\mathcal W_{\psi}^{\mathbb H}[f](a,{\bf y},\theta)&= f\otimes^{A_s} \psi\big({\bf x}\big)\\
&= a^{-1} \int_{\mathbb R^{2}} f({\bf x})\,e^{\frac{ja_2}{2b_2}(x_2^2-y_2^2)}\, \overline{\psi\big(r_{-\theta}a^{-1}({\bf x}-{\bf y})\big)}\,e^{\frac{ia_1}{2b_1}(x_1^2-y_1^2)} \,d{\bf x}\\
&=\left\langle f,\,\Psi_{a,{\bf y},\theta}^{\mathbb H}\right\rangle_{L^2\left(\mathbb R^{2},\,\mathbb H\right)}\tag{3.9}\\
\end{align*}
where $\Psi_{a,{\bf y},\theta}^{\mathbb H}$ is given by (3.5).

\parindent=0mm \vspace{.1in}
It is worth to note that the proposed quaternion linear canonical wavelet transform (3.9) boils down to existing quaternion wavelet transform as well as gives birth to some new quaternion wavelet transform which are not yet reported in the open literature:

\begin{itemize}

\item For the matrices $A_s =\left[\begin{array}{cc}0 & 1 \\-1 & 0 \\ \end{array}\right] $  for $s=1,2$,  the quaternion linear canonical wavelet transform (3.9) boils down to the quaternion wavelet transform given by \cite{Maw}
  \begin{align*}
  \mathcal W_{\psi}^{\mathbb H}[f](a,{\bf y},\theta)= a^{-1} \int_{\mathbb R^{2}} f({\bf x})  \overline{\psi\big(r_{-\theta}a^{-1}({\bf x}-{\bf y})\big)}\,d{\bf x}.\tag{3.10}
  \end{align*}

  \item For the matrices $A_s =\left[\begin{array}{cc}1 & b_s \\0 & 1 \\ \end{array}\right],\,b_s\neq 0 $, for $s=1,2$, we can obtain a new quaternion wavelet transform namely the quaternion Fresnel-Canonical wavelet transform given by
  \begin{align*}
  \mathcal W_{\psi}^{\mathbb H}[f](a,{\bf y},\theta)=a^{-1} \int_{\mathbb R^{2}} f({\bf x})\,e^{\frac{j}{2b_2}(x_2^2-y_2^2)} \,\overline{\psi\big(r_{-\theta}a^{-1}({\bf x}-{\bf y})\big)}\,e^{\frac{i}{2b_1}(x_1^2-y_1^2)} \,d{\bf x}.\tag{3.11}
  \end{align*}

  \item For the matrices $A_s =\left[\begin{array}{cc}\cos\alpha & \sin\alpha \\-\sin\alpha & \cos\alpha \\ \end{array}\right]$ for $s=1,2$,  the quaternion linear canonical wavelet transform (3.9) reduces to the quaternion fractional wavelet transform given by
  \begin{align*}
  \qquad\quad\mathcal W_{\psi}^{\mathbb H}[f](a,{\bf y},\theta)= a^{-1} \int_{\mathbb R^{2}} f({\bf x})\,e^{\frac{j\cot\alpha}{2}(x_2^2-y_2^2)} \, \overline{\psi\big(r_{-\theta}a^{-1}({\bf x}-{\bf y})\big)}\,e^{\frac{i\cot\alpha}{2}(x_1^2-y_1^2)}\,d{\bf x}.\tag{3.12}
  \end{align*}

\end{itemize}

\parindent=8mm \vspace{.1in}

For the demonstration of the quaternion linear canonical wavelet transform (3.9), we shall present an illustrative example.

\parindent=0mm \vspace{.1in}
{\bf Example 3.7.} Consider the two dimensional difference-of-Gaussian wavelets
\begin{align*}
\psi({\bf x})=\lambda^{-2} e^{\frac{-1}{2\lambda^2}(x_1^2+x_2^2)}\,-\,e^{\frac{-1}{2}(x_1^2+x_2^2)},\,0<\lambda<1.\tag{3.13}
\end{align*}
For $\theta=0,$ i.e. $r_{-\theta}=1,$ we have the corresponding quaternion linear canonical wavelet family of $\psi$ as

\begin{align*}
\Psi_{a,{\bf y},\theta}^{\mathbb H}({\bf x})&=a^{-1}\, e^{\frac{-ia_1}{2b_1}(x_1^2-y_1^2)} \psi\big(a^{-1}({\bf x}-{\bf y})\big)\,e^{\frac{-ja_2}{2b_2}(x_2^2-y_2^2)}\\
&=a^{-1} \,e^{\frac{-ia_1}{2b_1}(x_1^2-y_1^2)}  \left\{\lambda^{-2}e^{\frac{-\big((x_1-y_1)^2+(x_2-y_2)^2\big)}{2a^2\lambda^2}}\,-\,e^{\frac{-\big((x_1-y_1)^2+(x_2-y_2)^2\big)}{2a^2}}\right\}
e^{\frac{-ja_2}{2b_2}(x_2^2-y_2^2)}\\
\end{align*}

\parindent=0mm \vspace{.1in}
Then the quaternion linear canonical wavelet transform of the signal
\begin{align*}
f({\bf x})= e^{-\alpha_1x_1-\alpha_2x_2},\,\alpha_1,\alpha_2 \in \mathbb R
\end{align*}

with respect to above DOG wavelet $\psi$ (as defined in (3.13)) is given by
\begin{align*}
&\mathcal W_{\psi}^{\mathbb H}[f](a,{\bf y},\theta)\\
&\qquad= \int_{\mathbb R^{2}} f({\bf x})\,\overline{\Psi_{a,{\bf y},\theta}^{\mathbb H}({\bf x})}\,d{\bf x}\\
&\qquad= a^{-1}\int_{\mathbb R^{2}} e^{-\alpha_1x_1-\alpha_2x_2}\,e^{\frac{ja_2}{2b_2}(x_2^2-y_2^2)}\overline{   \left\{\lambda^{-2}e^{\frac{-\big((x_1-y_1)^2+(x_2-y_2)^2\big)}{2a^2\lambda^2}}\,-\,e^{\frac{-\big((x_1-y_1)^2+(x_2-y_2)^2\big)}{2a^2}}\right\}
}\,e^{\frac{ia_1}{2b_1}(x_1^2-y_1^2)}\,d{\bf x}\\
&\qquad= a^{-1}e^{\frac{-ja_2y_2^2}{2b_2}}\int_{\mathbb R^{2}} e^{-\alpha_1x_1-\alpha_2x_2}\,e^{\frac{ja_2x_2^2}{2b_2}}  \left\{\lambda^{-2}e^{\frac{-\big(x_1^2+y_1^2-2x_1y_1+x_2^2+y_2^2-2x_2y_2\big)}{2a^2\lambda^2}}\right.\\
&\qquad\qquad\qquad\qquad\qquad\qquad\quad\qquad\qquad\qquad\left.-\,e^{\frac{-\big(x_1^2+y_1^2-2x_1y_1)+(x_2^2+y_2^2-2x_2y_2\big)}{2a^2}}\right\}
\,e^{\frac{ia_1x_1^2}{2b_1}}\,d{\bf x} \,e^{\frac{-ia_1y_1^2}{2b_1}}\\
&\qquad= (\lambda^{2}a)^{-1}e^{\frac{-ja_2y_2^2}{2b_2}}\int_{\mathbb R^{2}} e^{-\alpha_2x_2}\,e^{\frac{ja_2x_2^2}{2b_2}}  e^{\frac{-\big(x_1^2+y_1^2-2x_1y_1+x_2^2+y_2^2-2x_2y_2\big)}{2a^2\lambda^2}}\,e^{-\alpha_1x_1}\,e^{\frac{ia_1x_1^2}{2b_1}}\,d{\bf x} \,e^{\frac{-ia_1y_1^2}{2b_1}}\\
&\qquad\qquad-a^{-1}e^{\frac{-ja_2y_2^2}{2b_2}}\int_{\mathbb R^{2}} e^{-\alpha_2x_2}\,e^{\frac{ja_2x_2^2}{2b_2}}\,e^{\frac{-\big(x_1^2+y_1^2-2x_1y_1+x_2^2+y_2^2-2x_2y_2\big)}{2a^2}}
\,e^{-\alpha_1x_1}\,e^{\frac{ia_1x_1^2}{2b_1}}\,d{\bf x} \,e^{\frac{-ia_1y_1^2}{2b_1}}\\
&\qquad= (\lambda^{2}a)^{-1}e^{\frac{-ja_2y_2^2}{2b_2}}\,e^{-\frac{y_1^2+y_2^2}{2a^2\lambda^2}}\int_{\mathbb R^{2}} e^{-\alpha_2x_2}\,e^{\frac{ja_2x_2^2}{2b_2}}  e^{\frac{-\big(x_1^2+\-2x_1y_1+x_2^2+-2x_2y_2\big)}{2a^2\lambda^2}}\,e^{-\alpha_1x_1}\,e^{\frac{ia_1x_1^2}{2b_1}}\,d{\bf x} \,e^{\frac{-ia_1y_1^2}{2b_1}}\\
&\qquad\qquad-a^{-1}e^{\frac{-ja_2y_2^2}{2b_2}}\,e^{-\frac{y_1^2+y_2^2}{2a^2}}\int_{\mathbb R^{2}} e^{-\alpha_2x_2}\,e^{\frac{ja_2x_2^2}{2b_2}}\,e^{\frac{-\big(x_1^2+y_1^2-2x_1y_1+x_2^2+y_2^2-2x_2y_2\big)}{2a^2}}
\,e^{-\alpha_1x_1}\,e^{\frac{ia_1x_1^2}{2b_1}}\,d{\bf x} \,e^{\frac{-ia_1y_1^2}{2b_1}}\\
&\qquad= \dfrac{1}{\lambda^{2}a}e^{\frac{-ja_2y_2^2}{2b_2}}\,e^{-\frac{y_1^2+y_2^2}{2a^2\lambda^2}}\int_{\mathbb R} \exp\left\{-x_2^2\left(\frac{1}{2a^2\lambda^2}-\frac{ja_2}{2b_2}\right)+x_2\left(-\alpha_2+\frac{y_2}{a^2\lambda^2}\right)\right\}dx_2\\
&\qquad\qquad\qquad\times\int_{\mathbb R}\exp\left\{-x_1^2\left(\frac{1}{2a^2\lambda^2}-\frac{ia_1}{2b_1}\right)+x_1\left(-\alpha_1+\frac{y_1}{a^2\lambda^2}\right)\right\}dx_1
\,e^{\frac{-ia_1y_1^2}{2b_1}}\\
&\qquad\qquad-\dfrac{1}{a}\,e^{\frac{-ja_2y_2^2}{2b_2}}\,e^{-\frac{y_1^2+y_2^2}{2a^2}}\int_{\mathbb R} \exp\left\{-x_2^2\left(\frac{1}{2a^2}-\frac{ja_2}{2b_2}\right)+x_2\left(-\alpha_2+\frac{y_2}{a^2}\right)\right\}dx_2\\
&\qquad\qquad\qquad\times\int_{\mathbb R}\exp\left\{-x_1^2\left(\frac{1}{2a^2}-\frac{ia_1}{2b_1}\right)+x_1\left(-\alpha_1+\frac{y_1}{a^2}\right)\right\}dx_1
\,e^{\frac{-ia_1y_1^2}{2b_1}}\\
\end{align*}
\begin{align*}
&\qquad= \dfrac{1}{\lambda^{2}a}\exp\left\{\frac{-ja_2y_2^2}{2b_2}-\frac{y_1^2+y_2^2}{2a^2\lambda^2}\right\}\sqrt{\dfrac{2\pi a^2\lambda^2b_2}{b_2-ja_2a^2\lambda^2}}\exp\left\{\frac{(y_2-a^2\alpha_2\lambda^2)^2}{a^4\lambda^4}\times\frac{2a^2\lambda^2b_2}{4(b_2-ja_2a^2\lambda^2)}\right\}\\
&\qquad\qquad\qquad\times \sqrt{\dfrac{2\pi a^2\lambda^2b_1}{b_1-ia_1a^2\lambda^2}}\exp\left\{\frac{(y_1-a^2\alpha_1\lambda^2)^2}{a^4\lambda^4}\times\frac{2a^2\lambda^2b_1}{4(b_1-ia_1a^2\lambda^2)}\right\}
\,e^{\frac{-ia_1y_1^2}{2b_1}}\\
&\qquad\qquad-\dfrac{1}{a}\exp\left\{\frac{-ja_2y_2^2}{2b_2}-\frac{y_1^2+y_2^2}{2a^2}\right\}\sqrt{\dfrac{2\pi a^2b_2}{b_2-ja_2a^2}}\exp\left\{\frac{(y_2-a^2\alpha_2)^2}{a^4}\times\frac{2a^2b_2}{4(b_2-ja_2a^2)}\right\}\\
&\qquad\qquad\qquad\times \sqrt{\dfrac{2\pi a^2b_1}{b_1-ia_1a^2}}\exp\left\{\frac{(y_1-a^2\alpha_1)^2}{a^4}\times\frac{2a^2b_1}{4(b_1-ia_1a^2)}\right\}
\,e^{\frac{-ia_1y_1^2}{2b_1}}\\
\end{align*}
\begin{align*}
&\qquad= \dfrac{2\pi a\sqrt{b_1b_2}}{\sqrt{b_2-ja_2a^2\lambda^2}} \exp\left\{\frac{-ja_2y_2^2}{2b_2}-\frac{y_1^2+y_2^2}{2a^2\lambda^2}+\frac{b_2(y_2-a^2\alpha_2\lambda^2)^2}{a^2\lambda^2(b_2-ja_2a^2\lambda^2)}\right\}
\sqrt{\dfrac{1}{b_1-ia_1a^2\lambda^2}}\\
&\qquad\qquad\qquad\times \exp\left\{\frac{b_1(y_1-a^2\alpha_1\lambda^2)^2}{a^2\lambda^2(b_1-ia_1a^2\lambda^2)}-\frac{ia_1y_1^2}{2b_1}\right\}-\dfrac{2\pi a\sqrt{b_1b_2}}{\sqrt{b_2-ja_2a^2}} \\
&\qquad\qquad\times \exp\left\{\frac{-ja_2y_2^2}{2b_2}-\frac{y_1^2+y_2^2}{2a^2}+\frac{b_2(y_2-a^2\alpha_2)^2}{a^2(b_2-ja_2a^2)}\right\} \sqrt{\dfrac{1}{b_1-ia_1a^2}}\exp\left\{\frac{b_1(y_1-a^2\alpha_1)^2}{a^2(b_1-ia_1a^2)}-\frac{ia_1y_1^2}{2b_1}\right\}\\
\end{align*}


\parindent=0mm \vspace{.0in}
In the following theorem, we assemble some of the basic properties of the proposed quaternion linear canonical wavelet transform (3.9).

\parindent=0mm \vspace{.1in}
{\bf Theorem 3.8.} {\it Let  $\psi,\phi\in L^2(\mathbb R^2,\mathbb H)$ be two admissible wavelets, then for every $f,g\in L^2(\mathbb R^2,\mathbb H)$ linear canonical wavelet transform (3.9) satisfies the following properties:}

\parindent=0mm \vspace{.1in}
(i)~~${\text {\it Linearity:}}\qquad{\mathcal W}_{\psi}^{\mathbb H}\Big[ \alpha_1 f+\alpha_2 g\Big](a,{\bf y},\theta)=\alpha_1{\mathcal W}_{\psi}^{\mathbb H}\big[ f\big](a,{\bf y},\theta)+\alpha_2 {\mathcal W}_{\psi}^{\mathbb H}\big[g\big](a,{\bf y},\theta),$\\

$\quad\text{where } \alpha_1, \alpha_2 \text{are quaternion constants in }{\mathbb H}$

\parindent=0mm \vspace{.1in}
(ii)~~${\text {\it Anti-linearity:}}\qquad{\mathcal W}_{\alpha_1 \psi +\alpha_2 \phi}^{\mathbb H}\big[ f\big](a,{\bf y},\theta)={\mathcal W}_{\psi}^{\mathbb H}\big[ f\big](a,{\bf y},\theta)\overline{\alpha_1}+ {\mathcal W}_{\phi}^{\mathbb H}\big[g\big](a,{\bf y},\theta)\overline{\alpha_2},$

\parindent=0mm \vspace{.1in}
(iii)~${\text {\it Translation:}}~~\text{For any constant {\bf k}}\in\mathbb R^2,$\\
$${\mathcal W}_{\psi}^{\mathbb H}\big[ f({\bf x}-{\bf k})\big](a,{\bf y},\theta)= e^{\frac{ia_1k_1^2}{b_1}}\,{\mathcal W}_{\psi}^{\mathbb H}\Big[ e^{\frac{ia_1k_1x_1}{b_1}}\,f({\bf x}) \,e^{\frac{ja_2k_2x_2}{b_2}}\Big](a,{\bf y}-{\bf k},\theta)\,e^{\frac{ja_2k_2^2}{b_2}}$$

\parindent=0mm \vspace{.0in}
(iv)\,${\text {\it Scaling:}} \,\,\text{For non-zero constant }\lambda,$\\
$${\mathcal W}_{\psi}^{\mathbb H}\big[ f(\lambda\,{\bf x})\big](a,{\bf y},\theta)=\dfrac{1}{\lambda}{\mathcal W}_{\psi}^{\mathbb H}\big[ f(\lambda\,{\bf x})\big](a\lambda,{\bf y}\lambda,\theta),$$
where matrices parameter in R.H.S is $A^{\prime}_s=[a_s^{\prime},b_s^{\prime},c_s,d_s]$ is related to matrices in L.H.S i.e. $A_s=[a_s,b_s,c_s,d_s]$ by a relation $\dfrac{a_s^{\prime}}{b_s^{\prime}}=\dfrac{a_s}{b_s\lambda^2}$, for $s=1,2.$

\parindent=0mm \vspace{.0in}
(v)\,${\text {\it Parity:}}\quad{\mathcal W}_{P\psi}^{\mathbb H}\big[ Pf({\bf x})\big](a,{\bf y},\theta)={\mathcal W}_{\psi}^{\mathbb H}\big[ f({\bf x})\big](a,-{\bf y},\theta),\,\text{where}\, Pf({\bf x})=f(-{\bf x}).$

\parindent=0mm \vspace{.0in}
(vi)\,${\text {\it Dilation in $\psi$:}}\quad{\mathcal W}_{D_c\psi}^{\mathbb H}\big[ f({\bf x})\big](a,{\bf y},\theta)={\mathcal W}_{\psi}^{\mathbb H}\big[ f({\bf x})\big](ac,{\bf y},\theta),\,\text{where}\, D_c\psi({\bf x})=\frac{1}{c}\psi\left(\frac{{\bf x}}{c}\right).$

\parindent=0mm \vspace{.1in}
{\it Proof.} For brevity proof of these properties is omitted.

\parindent=0mm \vspace{.1in}
In our next theorem, we will show that the quaternion linear canonical wavelet transform sets up an isometry from $L^2(\mathbb R^2,\mathbb H)$ to $L^2(\mathbb R^{+}\times\mathbb R^2\times\text{SO(2)},\mathbb H)$.

\parindent=0mm \vspace{.1in}
{\bf Theorem 3.9(Parseval's Formula).} {\it Suppose that $\psi\in L^2(\mathbb R^2,\mathbb H)$ be an admissible wavelets, then for every $f,g\in L^2(\mathbb R^2,\mathbb H),$ we have }

\begin{align*}
\left\langle\mathcal W_{\psi}^{\mathbb H}\big[ f\big],\,\mathcal W_{\psi}^{\mathbb H}\big[ g\big]\right\rangle_{L^2(\mathcal G,\mathbb H)}\,=\, C_{\psi} \langle f,\,g\rangle_{L^2(\mathbb R^2,\mathbb H)}\,\tag{3.13}
\end{align*}

where $\mathcal G \,=\mathbb R^{+}\times\mathbb R^2\times\text{SO(2)}$ is similitude group and $C_{\psi}$ is admissiblity given by (3.8).

\parindent=0mm \vspace{.1in}
{\it Proof.} Invoking Parseval formula (2.4) for the two-sided QLCT and implementing Lemma 3.4, we have

\begin{align*}
&\left\langle\mathcal W_{\psi}^{\mathbb H}\big[ f\big],\,\mathcal W_{\psi}^{\mathbb H}\big[ g\big]\right\rangle_{L^2(\mathcal G,\mathbb H)}\\
&\quad= \int_{\mathcal G } \mathcal W_{\psi}^{\mathbb H}\big[ f\big](a,{\bf y},\theta)\,\overline{\mathcal W_{\psi}^{\mathbb H}\big[ g\big](a,{\bf y},\theta)}\,\dfrac{da d{\bf y} d\theta}{a^3} \\
&\quad= \int_{\mathcal G } \left\langle f,\,\Psi_{a,{\bf y},\theta}^{\mathbb H}\right\rangle_{L^2\left(\mathbb R^{2},\,\mathbb H\right)}\,\overline{\left\langle g,\,\Psi_{a,{\bf y},\theta}^{\mathbb H}\right\rangle}_{L^2\left(\mathbb R^{2},\,\mathbb H\right)}\,\dfrac{da d{\bf y} d\theta}{a^3} \\
&\quad= \int_{\mathcal G } \left\langle \mathcal L^{\mathbb H}_{A_1,A_2}\big[f\big],\,\mathcal L^{\mathbb H}_{A_1,A_2}\Big[\Psi_{a,{\bf y},\theta}^{\mathbb H}({\bf x})\Big]\right\rangle_{L^2\left(\mathbb R^{2},\,\mathbb H\right)}\,\overline{\left\langle \mathcal L^{\mathbb H}_{A_1,A_2}\big[g\big],\,\mathcal L^{\mathbb H}_{A_1,A_2}\Big[\Psi_{a,{\bf y},\theta}^{\mathbb H}({\bf x})\Big]\right\rangle}_{L^2\left(\mathbb R^{2},\,\mathbb H\right)}\,\dfrac{da d{\bf y} d\theta}{a^3} \\
&\quad= \int_{\mathcal G } \int_{\mathbb R^2}\mathcal L^{\mathbb H}_{A_1,A_2}\big[f\big]({\bf w})\overline{\mathcal L^{\mathbb H}_{A_1,A_2}\Big[\Psi_{a,{\bf y},\theta}^{\mathbb H}({\bf x})\Big]}({\bf w})d{\bf w}\,\int_{\mathbb R^2}\mathcal L^{\mathbb H}_{A_1,A_2}\Big[\Psi_{a,{\bf y},\theta}^{\mathbb H}({\bf x})\Big]({\bf w}^{\prime})\overline{\mathcal L^{\mathbb H}_{A_1,A_2}\big[g\big]({\bf w}^{\prime})}d{\bf w}^{\prime}\dfrac{da d{\bf y} d\theta}{a^3} \\
&\quad= \dfrac{a^2}{(2\pi)^2 b_1b_2}\int_{\mathcal G } \int_{\mathbb R^2}\int_{\mathbb R^2}\mathcal L^{\mathbb H}_{A_1,A_2}\big[f\big]({\bf w})\\
&\quad\times\overline{ \exp\left\{\frac{i}{2b_1}\big[d_1w_1^2+y_1^2a_1-2y_1w_1-\frac{\pi b_1}{2}\big]\right\} \mathcal F_q[\psi]\left(r_{-\theta}\,\dfrac{a{\bf w}}{{\bf b}}\right)\exp\left\{\frac{j}{2b_2}\big[d_2w_2^2+y_2^2a_2-2y_2w_2-\frac{\pi b_2}{2}\big]\right\}}\\
&\quad\quad\times \exp\left\{\frac{i}{2b_1}\big[d_1w^{\prime^2}_1+y_1^2a_1-2y_1 w^{\prime}_1-\frac{\pi b_1}{2}\big]\right\} \mathcal F_q[\psi]\left(r_{-\theta}\,\dfrac{a{\bf w^{\prime}}}{{\bf b}}\right)\\
&\quad\quad\times\exp\left\{\frac{j}{2b_2}\big[d_2 w^{\prime^2}_2+y_2^2a_2-2y_2w^{\prime}_2-\frac{\pi b_2}{2}\big]\right\}\,\overline{\mathcal L^{\mathbb H}_{A_1,A_2}\big[g\big]({\bf w}^{\prime})},d{\bf w} d{\bf w}^{\prime}\,\dfrac{da d{\bf y} d\theta}{a^3} \\
&\quad= \dfrac{a^2}{(2\pi)^2 b_1b_2}\int_{\mathcal G } \int_{\mathbb R^2}\int_{\mathbb R^2}\mathcal L^{\mathbb H}_{A_1,A_2}\big[f\big]({\bf w})\exp\left\{\frac{-j}{2b_2}\big[d_2w_2^2+y_2^2a_2-2y_2w_2-\frac{\pi b_2}{2}\big]\right\}\\
&\quad\quad\times\overline{ \mathcal F_q[\psi]\left(r_{-\theta}\,\dfrac{a{\bf w}}{{\bf b}}\right)}\exp\left\{\frac{-i}{2b_1}\big[d_1w_1^2+y_1^2a_1-2y_1w_1-\frac{\pi b_1}{2}-d_1 w^{\prime^2}_1-y_1^2a_1+2y_1 w^{\prime}_1+\frac{\pi b_1}{2}\big]\right\}\\
&\quad\quad\times  \mathcal F_q[\psi]\left(r_{-\theta}\,\dfrac{a{\bf w^{\prime}}}{{\bf b}}\right)\exp\left\{\frac{j}{2b_2}\big[d_2 w^{\prime^2}_2+y_2^2a_2-2y_2 w^{\prime}_2-\frac{\pi b_2}{2}\big]\right\}\,\overline{\mathcal L^{\mathbb H}_{A_1,A_2}\big[g\big]({\bf w}^{\prime})}\,d{\bf w} d{\bf w}^{\prime}\,\dfrac{da d{\bf y} d\theta}{a^3} \\
&\quad= \dfrac{a^2}{(2\pi)^2 b_1b_2}\int_{\mathcal G } \int_{\mathbb R^2}\int_{\mathbb R^2}\mathcal L^{\mathbb H}_{A_1,A_2}\big[f\big]({\bf w})\exp\left\{\frac{-j}{2b_2}\big[d_2w_2^2+y_2^2a_2-2y_2w_2-\frac{\pi b_2}{2}\big]\right\}\\
&\quad\quad\times\overline{ \mathcal F_q[\psi]\left(r_{-\theta}\,\dfrac{a{\bf w}}{{\bf b}}\right)}\exp\left\{\frac{-i}{2b_1}\big[d_1(w_1^2-w^{\prime^2}_1)-2y_1(w_1-w^{\prime}_1)\big]\right\}\\
&\quad\quad\times  \mathcal F_q[\psi]\left(r_{-\theta}\,\dfrac{a{\bf w^{\prime}}}{{\bf b}}\right)\exp\left\{\frac{j}{2b_2}\big[d_2 w^{\prime^2}_2+y_2^2a_2-2y_2 w^{\prime}_2-\frac{\pi b_2}{2}\big]\right\}\,\overline{\mathcal L^{\mathbb H}_{A_1,A_2}\big[g\big]({\bf w}^{\prime})}\,d{\bf w} d{\bf w}^{\prime}\,\dfrac{da d{\bf y} d\theta}{a^3} \\
&\quad= \dfrac{a^2}{(2\pi)^2 b_1b_2}\int_{\mathbb R^{+}\times\mathbb R^2\times\text{SO(2)} } \int_{\mathbb R^2}\int_{\mathbb R}\mathcal L^{\mathbb H}_{A_1,A_2}\big[f\big]({\bf w})\exp\left\{\frac{-j}{2b_2}\big[d_2w_2^2+y_2^2a_2-2y_2w_2-\frac{\pi b_2}{2}\big]\right\}\\
&\quad\quad\times\overline{ \mathcal F_q[\psi]\left(r_{-\theta}\,\dfrac{a{\bf w}}{{\bf b}}\right)}\exp\left\{\frac{-id_1}{2b_1}(w_1^2- w^{\prime^2}_1)\right\}
\int_{\mathbb R}\exp\left\{\frac{iy_1}{b_1}(w_1- w^{\prime}_1)\right\}\,dy_1\\
&\quad\quad\times  \mathcal F_q[\psi]\left(r_{-\theta}\,\dfrac{a{\bf w^{\prime}}}{{\bf b}}\right)\exp\left\{\frac{j}{2b_2}\big[d_2 w^{\prime^2}_2+y_2^2a_2-2y_2  w^{\prime}_2-\frac{\pi b_2}{2}\big]\right\}\,\,\overline{\mathcal L^{\mathbb H}_{A_1,A_2}\big[g\big]({\bf w}^{\prime})}d{\bf w} d{\bf w}^{\prime}\,\dfrac{da dy_2 d\theta}{a^3} \\
&\quad= \dfrac{a^2}{(2\pi)^2 b_1b_2}\int_{\mathbb R^{+}\times\mathbb R^2\times\text{SO(2)} } \int_{\mathbb R^2}\int_{\mathbb R}\mathcal L^{\mathbb H}_{A_1,A_2}\big[f\big]({\bf w})\exp\left\{\frac{-j}{2b_2}\big[d_2w_2^2+y_2^2a_2-2y_2w_2-\frac{\pi b_2}{2}\big]\right\}\\
&\quad\qquad\times\overline{ \mathcal F_q[\psi]\left(r_{-\theta}\,\dfrac{a{\bf w}}{{\bf b}}\right)}\exp\left\{\frac{-id_1}{2b_1}(w_1^2- w^{\prime^2}_1)\right\}
\delta(w_1- w^{\prime}_1)\mathcal F_q[\psi]\left(r_{-\theta}\,\dfrac{a{\bf w^{\prime}}}{{\bf b}}\right)\\
&\quad\qquad\times  \exp\left\{\frac{j}{2b_2}\big[d_2 w^{\prime^2}_2+y_2^2a_2-2y_2 w^{\prime}_2-\frac{\pi b_2}{2}\big]\right\}\,\overline{\mathcal L^{\mathbb H}_{A_1,A_2}\big[g\big]({\bf w}^{\prime})}\, d w^{\prime}_2 dw_1dw_2\,\dfrac{da dy_2 d\theta}{a^3} \\
&\quad= \dfrac{a^2}{(2\pi)^2 b_1b_2}\int_{\mathbb R^{+}\times\mathbb R^2\times\text{SO(2)} } \int_{\mathbb R^2}\mathcal L^{\mathbb H}_{A_1,A_2}\big[f\big]({\bf w})\exp\left\{\frac{-j}{2b_2}\big[d_2w_2^2+y_2^2a_2-2y_2w_2-\frac{\pi b_2}{2}\big]\right\}\\
&\quad\qquad\times\overline{ \mathcal F_q[\psi]\left(r_{-\theta}\,\dfrac{a{\bf w}}{{\bf b}}\right)}\mathcal F_q[\psi]\left(r_{-\theta}\,\dfrac{a{\bf w^{\prime}}}{{\bf b}}\right)\exp\left\{\frac{j}{2b_2}\big[d_2 w^{\prime^2}_2+y_2^2a_2-2y_2 w^{\prime}_2-\frac{\pi b_2}{2}\big]\right\}\\
&\quad\qquad\times\overline{\mathcal L^{\mathbb H}_{A_1,A_2}\big[g\big]({\bf w}^{\prime})}\, d w^{\prime}_2 dw_1dw_2\,\dfrac{da dy_2 d\theta}{a^3} \\
\end{align*}

\begin{align*}
&\quad= \dfrac{a^2}{(2\pi)^2 b_1b_2}\int_{\mathbb R^{+}\times\mathbb R^2\times\text{SO(2)} } \int_{\mathbb R^2}\mathcal L^{\mathbb H}_{A_1,A_2}\big[f\big]({\bf w})\overline{\exp\left\{\frac{j}{2b_2}\big[d_2w_2^2+y_2^2a_2-2y_2w_2-\frac{\pi b_2}{2}\big]\right\}}\\
&\quad\qquad\times\overline{ \mathcal F_q[\psi]\left(r_{-\theta}\,\dfrac{a{\bf w}}{{\bf b}}\right)}\mathcal F_q[\psi]\left(r_{-\theta}\,\dfrac{a{\bf w^{\prime}}}{{\bf b}}\right)\exp\left\{\frac{j}{2b_2}\big[d_2 w^{\prime^2}_2+y_2^2a_2-2y_2 w^{\prime}_2-\frac{\pi b_2}{2}\big]\right\}\\
&\quad\qquad\times\overline{\mathcal L^{\mathbb H}_{A_1,A_2}\big[g\big]({\bf w}^{\prime})}\, d w^{\prime}_2 dw_1dw_2\,\dfrac{da dy_2 d\theta}{a^3} \\
&\quad= \dfrac{a^2}{(2\pi)^2 b_1b_2}\int_{\mathbb R^{+}\times\mathbb R^2\times\text{SO(2)} } \int_{\mathbb R^2}\mathcal L^{\mathbb H}_{A_1,A_2}\big[f\big]({\bf w})\overline{\exp\left\{\frac{j}{2b_2}\big[d_2w_2^2+y_2^2a_2-2y_2w_2-\frac{\pi b_2}{2}\big]\right\}}\\
&\quad\qquad\times\overline{ \mathcal F_q[\psi]\left(r_{-\theta}\,\dfrac{a{\bf w}}{{\bf b}}\right)}\,\,\overline{\overline{\mathcal F_q[\psi]\left(r_{-\theta}\,\dfrac{a{\bf w^{\prime}}}{{\bf b}}\right)}}\overline{\exp\left\{\frac{-j}{2b_2}\big[d_2 w^{\prime^2}_2+y_2^2a_2-2y_2 w^{\prime}_2-\frac{\pi b_2}{2}\big]\right\}}\\
&\quad\qquad\times\overline{\mathcal L^{\mathbb H}_{A_1,A_2}\big[g\big]({\bf w}^{\prime})}\, d w^{\prime}_2 dw_1dw_2\,\dfrac{da dy_2 d\theta}{a^3} \\
&\quad= \dfrac{a^2}{(2\pi)^2 b_1b_2}\int_{\mathbb R^{+}\times\mathbb R^2\times\text{SO(2)} } \int_{\mathbb R^2}\mathcal L^{\mathbb H}_{A_1,A_2}\big[f\big]({\bf w})\overline{\mathcal F_q[\psi]\left(r_{-\theta}\,\dfrac{a{\bf w}}{{\bf b}}\right)\exp\left\{\frac{j}{2b_2}\big[d_2w_2^2+y_2^2a_2-2y_2w_2-\frac{\pi b_2}{2}\big]\right\}}\\
&\quad\quad\times\,\,\overline{\exp\left\{\frac{-j}{2b_2}\big[d_2{\bf w}^{\prime^2}_2+y_2^2a_2-2y_2{\bf w}^{\prime}_2-\frac{\pi b_2}{2}\big]\right\}\overline{\mathcal F_q[\psi]\left(r_{-\theta}\,\dfrac{a{\bf w^{\prime}}}{{\bf b}}\right)}}\,\overline{\mathcal L^{\mathbb H}_{A_1,A_2}\big[g\big]({\bf w}^{\prime})}\, d w^{\prime}_2 dw_1dw_2\,\dfrac{da dy_2 d\theta}{a^3} \\
&\quad= \dfrac{a^2}{(2\pi)^2 b_1b_2}\int_{\mathbb R^{+}\times\mathbb R^2\times\text{SO(2)} } \int_{\mathbb R}\mathcal L^{\mathbb H}_{A_1,A_2}\big[f\big]({\bf w})\\
&\quad\qquad\times\overline{\mathcal F_q[\psi]\left(r_{-\theta}\,\dfrac{a{\bf w}}{{\bf b}}\right)\exp\left\{\frac{jd_2}{2b_2}(w_2^2-{\bf w}^{\prime^2}_2)\right\}
\int_{\mathbb R}\exp\left\{\frac{-jy_2}{b_2}(w_2-{\bf w}^{\prime}_2)\right\}dy_2 \overline{\mathcal F_q[\psi]\left(r_{-\theta}\,\dfrac{a{\bf w^{\prime}}}{{\bf b}}\right)}}\\
&\quad\qquad\times\overline{\mathcal L^{\mathbb H}_{A_1,A_2}\big[g\big]({\bf w}^{\prime})}\, d w^{\prime}_2 dw_1dw_2\,\dfrac{da d\theta}{a^3} \\
&\quad= \dfrac{a^2}{(2\pi)^2 b_1b_2}\int_{\mathbb R^{+}\times\mathbb R^2\times\text{SO(2)} } \int_{\mathbb R}\mathcal L^{\mathbb H}_{A_1,A_2}\big[f\big]({\bf w})\\
&\quad\qquad\times\overline{\mathcal F_q[\psi]\left(r_{-\theta}\,\dfrac{a{\bf w}}{{\bf b}}\right)\exp\left\{\frac{jd_2}{2b_2}(w_2^2-{\bf w}^{\prime^2}_2)\right\}
\,\delta({\bf w}^{\prime}_2-w_2) \overline{\mathcal F_q[\psi]\left(r_{-\theta}\,\dfrac{a{\bf w^{\prime}}}{{\bf b}}\right)}}\\
&\quad\qquad\times\overline{\mathcal L^{\mathbb H}_{A_1,A_2}\big[g\big]({\bf w}^{\prime})} \, d w^{\prime}_2 dw_1dw_2\,\dfrac{da d\theta}{a^3} \\
&\quad= \dfrac{a^2}{(2\pi)^2 b_1b_2}\int_{\mathbb R^{+}\times\mathbb R^2\times\text{SO(2)} } \mathcal L^{\mathbb H}_{A_1,A_2}\big[f\big]({\bf w})\left|\mathcal F_q[\psi]\left(r_{-\theta}\,\dfrac{a{\bf w}}{{\bf b}}\right)\right|^2\,\overline{\mathcal L^{\mathbb H}_{A_1,A_2}\big[g\big]({\bf w})} \,\,dw_1dw_2 \,\dfrac{da d\theta}{a^3} \\
&\quad= \dfrac{a^2}{(2\pi)^2 b_1b_2}\int_{\mathbb R^{+}\times\mathbb R^2\times\text{SO(2)} } \mathcal L^{\mathbb H}_{A_1,A_2}\big[f\big]({\bf w})\left|\frac{2\pi \sqrt{b_1b_2}}{a}\mathcal L^{\mathbb H}_{A_1,A_2}\Big[\Psi_{a,{\bf y},\theta}^{\mathbb H}({\bf x})\Big]({\bf w})\right|^2\,\overline{\mathcal L^{\mathbb H}_{A_1,A_2}\big[g\big]({\bf w})} \, dw_1dw_2\,\dfrac{da d\theta}{a^3} \\
&\quad= \int_{\mathbb R^{+}\times\mathbb R^2\times\text{SO(2)} } \mathcal L^{\mathbb H}_{A_1,A_2}\big[f\big]({\bf w})\left|\mathcal L^{\mathbb H}_{A_1,A_2}\Big[\Psi_{a,{\bf y},\theta}^{\mathbb H}({\bf x})\Big]({\bf w})\right|^2\,\overline{\mathcal L^{\mathbb H}_{A_1,A_2}\big[g\big]({\bf w})} \,\,d{\bf w} \,\dfrac{da d\theta}{a^3} \\
&\quad= \int_{\mathbb R^2 } \mathcal L^{\mathbb H}_{A_1,A_2}\big[f\big]({\bf w})\int_{\mathbb R^{+}\times\text{SO(2)}}\left|\mathcal L^{\mathbb H}_{A_1,A_2}\Big[\Psi_{a,{\bf y},\theta}^{\mathbb H}({\bf x})\Big]({\bf w})\right|^2\,\dfrac{da  d\theta}{a^3} \,\overline{\mathcal L^{\mathbb H}_{A_1,A_2}\big[g\big]({\bf w})} \,d{\bf w} \\
&\quad=\, C_{\psi} \left\langle \mathcal L^{\mathbb H}_{A_1,A_2}\big[f\big],\,\mathcal L^{\mathbb H}_{A_1,A_2}\big[g\big]\right\rangle_{L^2(\mathbb R^2,\mathbb H)}\\
&\quad=\, C_{\psi} \langle f,\,g\rangle_{L^2(\mathbb R^2,\mathbb H)}
\end{align*}

where $C_\psi$ is given by (3.8).

\parindent=0mm \vspace{.1in}
This completes the proof of the theorem.\quad\fbox

\parindent=0mm \vspace{.1in}

{\bf Corollary 3.10(Energy Conservation).} {\it For $f=g$, we have the following identity:}
\begin{align*}
\int_{\mathbb R^{+}\times\mathbb R^2\times\text{SO(2)} }\Big\|\mathcal W_{\psi}^{\mathbb H}\big[ f\big](a,{\bf y},\theta)\Big\|_{\mathbb H}^{2}\dfrac{da\,d{\bf y}\,d\theta}{a^{3}}=C_{\psi}\, \big\|f\big\|^{2}_{L^2(\mathbb R^2,\mathbb H)}.\tag{3.14}
\end{align*}

{\it Remark:} We see that, except the factor $C_{\psi}$, the quaternion linear canonical wavelet transform sets up an isometry from $L^2(\mathbb R^2,\mathbb H)$ to $L^2(\mathbb R^{+}\times\mathbb R^2\times\text{SO(2)},\mathbb H)$.

\parindent=8mm \vspace{.1in}
The next theorem guarantees the reconstruction of the input quaternion  signal from the corresponding quaternion linear canonical wavelet transform.

\parindent=0mm \vspace{.1in}
{\bf Theorem 3.11(Inversion Formula).} {\it Suppose that $\psi\in L^2(\mathbb R^2,\mathbb H)$ be an admissible wavelets, then any quaternion  signal $f\in L^2(\mathbb R^2,\mathbb H)$ can be reconstructed from the quaternion linear canonical wavelet transform $\mathcal W_{\psi}^{\mathbb H}\big[ f\big](a,{\bf y},\theta)$ via the following formula:}
\begin{align*}
f({\bf x})=\dfrac{1}{C_{\Psi}}\int_{\mathbb R^{+}\times\mathbb R^2\times\text{SO(2)}}\mathcal W_{\psi}^{\mathbb H}\big[ f\big](a,{\bf y},\theta)\,\Psi_{a,{\bf y},\theta}^{\mathbb H}({\bf x})\,\dfrac{da\,d{\bf y}\,d\theta}{a^{3}},\,a.e.\,\tag{3.15}
\end{align*}

\parindent=0mm \vspace{.1in}
{\it Proof.} For arbitrary $g\in L^2(\mathbb R^2,\mathbb H),$ implication of Theorem 3.9 yields

\begin{align*}
 C_{\psi} \langle f,\,g\rangle_{L^2(\mathbb R^2,\mathbb H)} &=\int_{\mathcal G } \mathcal W_{\psi}^{\mathbb H}\big[ f\big](a,{\bf y},\theta)\,\overline{\mathcal W_{\psi}^{\mathbb H}\big[ g\big](a,{\bf y},\theta)}\,\dfrac{da d{\bf y} d\theta}{a^3}\\
 &=\int_{\mathcal G } \mathcal W_{\psi}^{\mathbb H}\big[ f\big](a,{\bf y},\theta)\,\overline{\left\langle g,\,\Psi_{a,{\bf y},\theta}^{\mathbb H}\right\rangle}_{L^2\left(\mathbb R^{2},\,\mathbb H\right)}\,\dfrac{da d{\bf y} d\theta}{a^3}\\
&=\int_{\mathcal G } \mathcal W_{\psi}^{\mathbb H}\big[ f\big](a,{\bf y},\theta)\,\left\langle \Psi_{a,{\bf y},\theta}^{\mathbb H},\,g\right\rangle_{L^2\left(\mathbb R^{2},\,\mathbb H\right)}\,\dfrac{da d{\bf y} d\theta}{a^3}\\
&=\int_{\mathcal G } \mathcal W_{\psi}^{\mathbb H}\big[ f\big](a,{\bf y},\theta)\,\int_{\mathbb R^2} \Psi_{a,{\bf y},\theta}^{\mathbb H},\overline{g({\bf x})}\,d{\bf x}\,\dfrac{da d{\bf y} d\theta}{a^3}\\
&=\int_{\mathbb R^2} \int_{\mathcal G } \mathcal W_{\psi}^{\mathbb H}\big[ f\big](a,{\bf y},\theta)\,\Psi_{a,{\bf y},\theta}^{\mathbb H}\,\dfrac{da d{\bf y} d\theta}{a^3}\overline{g({\bf x})}\,d{\bf x}\\
&=\left\langle \int_{\mathcal G } \mathcal W_{\psi}^{\mathbb H}\big[ f\big](a,{\bf y},\theta)\,\Psi_{a,{\bf y},\theta}^{\mathbb H}\,\dfrac{da d{\bf y} d\theta}{a^3},\,\,g\right\rangle_{L^2\left(\mathbb R^{2},\,\mathbb H\right)}
\end{align*}
where we applied Fobini's theorem in getting second last equality. Therefore, we have
\begin{align*}
C_{\Psi}\,f({\bf x})=\int_{\mathbb R^{+}\times\mathbb R^2\times\text{SO(2)}}\mathcal W_{\psi}^{\mathbb H}\big[ f\big](a,{\bf y},\theta)\,\Psi_{a,{\bf y},\theta}^{\mathbb H}({\bf x})\,\dfrac{da\,d{\bf y}\,d\theta}{a^{3}},\,a.e.
\end{align*}

This completes the proof of theorem. \quad\fbox

\parindent=0mm \vspace{.1in}

{\bf Theorem 3.12(Characterization of range).} For an admissible wavelet $\psi\in L^2(\mathbb R^2,\mathbb H)$, the range of the quaternion linear canonical wavelet transform $\mathcal W_{\psi}^{\mathbb H}$ is a reproducing kernel in $L^2\left(\mathbb R^{+}\times\mathbb R^2\times\text{SO(2)},\,\mathbb H\right)$ where kernel is given by

\begin{align*}
K_{\psi}(a,{\bf y},\theta; a^{\prime},{\bf y}^{\prime},\theta^{\prime} )=C_{\psi}^{-1}\left\langle\Psi_{a,{\bf y},\theta}^{\mathbb H},\,\,\Psi_{a^{\prime},{\bf y}^{\prime},\theta^{\prime}}^{\mathbb H} \right\rangle_{L^2\left(\mathbb R^{2},\,\mathbb H\right)}\,\tag{3.16}
\end{align*}

Moreover,
\begin{align*}
\left|K_{\psi}(a,{\bf y},\theta; a^{\prime},{\bf y}^{\prime},\theta^{\prime} )\right|\,\leq\,C_{\psi}^{-1}\big\|\psi\big\|_{L^2\left(\mathbb R^{2},\,\mathbb H\right)},\, \text{if}\, C_{\psi}>0.\,\tag{3.17}
\end{align*}

\parindent=0mm \vspace{.1in}
{\it Proof.} Invoking the inversion formula (3.15) in the definition of quaternion linear canonical wavelet transform (3.8),we have

\begin{align*}
\mathcal W_{\psi}^{\mathbb H}\big[f \big](a^{\prime},{\bf y}^{\prime},\theta^{\prime})&=\int_{\mathbb R^2} f({\bf x})\overline{\,\Psi_{a^{\prime},{\bf y}^{\prime},\theta^{\prime}}^{\mathbb H}({\bf x})}\,d{\bf x}\\
&=\int_{\mathbb R^2} \dfrac{1}{C_{\psi}}\int_{\mathbb R^{+}\times\mathbb R^2\times\text{SO(2)}}\mathcal W_{\psi}^{\mathbb H}\big[ f\big](a,{\bf y},\theta)\,\Psi_{a,{\bf y},\theta}^{\mathbb H}({\bf x})\,\dfrac{da\,d{\bf y}\,d\theta}{a^{3}}\overline{\,\Psi_{a^{\prime},{\bf y}^{\prime},\theta^{\prime}}^{\mathbb H}({\bf x})}\,d{\bf x}\\
&=\int_{\mathbb R^{+}\times\mathbb R^2\times\text{SO(2)}}\int_{\mathbb R^2} \mathcal W_{\psi}^{\mathbb H}\big[ f\big](a,{\bf y},\theta)\dfrac{1}{C_{\psi}}\,\Psi_{a,{\bf y},\theta}^{\mathbb H}({\bf x})\overline{\,\Psi_{a^{\prime},{\bf y}^{\prime},\theta^{\prime}}^{\mathbb H}({\bf x})}\,d{\bf x}\,\dfrac{da\,d{\bf y}\,d\theta}{a^{3}}\\
&=\int_{\mathbb R^{+}\times\mathbb R^2\times\text{SO(2)}} \mathcal W_{\psi}^{\mathbb H}\big[ f\big](a,{\bf y},\theta)\dfrac{1}{C_{\psi}}\,\int_{\mathbb R^2}\Psi_{a,{\bf y},\theta}^{\mathbb H}({\bf x})\overline{\,\Psi_{a^{\prime},{\bf y}^{\prime},\theta^{\prime}}^{\mathbb H}({\bf x})}\,d{\bf x}\,\dfrac{da\,d{\bf y}\,d\theta}{a^{3}}\\
&=\int_{\mathbb R^{+}\times\mathbb R^2\times\text{SO(2)}} \mathcal W_{\psi}^{\mathbb H}\big[ f\big](a,{\bf y},\theta)C_{\psi}^{-1}\left\langle\Psi_{a,{\bf y},\theta}^{\mathbb H},\,\,\Psi_{a^{\prime},{\bf y}^{\prime},\theta^{\prime}}^{\mathbb H} \right\rangle_{L^2\left(\mathbb R^{2},\,\mathbb H\right)}\,\dfrac{da\,d{\bf y}\,d\theta}{a^{3}}\\
&=\int_{\mathbb R^{+}\times\mathbb R^2\times\text{SO(2)}} \mathcal W_{\psi}^{\mathbb H}\big[ f\big](a,{\bf y},\theta)K_{\psi}(a,{\bf y},\theta; a^{\prime},{\bf y}^{\prime},\theta^{\prime} )\,\dfrac{da\,d{\bf y}\,d\theta}{a^{3}}\\
\end{align*}
Which completes the proof of first assertion.

\parindent=0mm \vspace{.1in}
Also, from the Definition 3.3, we have
\begin{align*}
\left|K_{\psi}(a,{\bf y},\theta; a^{\prime},{\bf y}^{\prime},\theta^{\prime} )\right|\,&\leq\,\dfrac{1}{\left|C_{\psi}\right|}\int_{\mathbb R^2} \left|\Psi_{a,{\bf y},\theta}^{\mathbb H}({\bf x})\,\overline{\Psi_{a^{\prime},{\bf y}^{\prime},\theta^{\prime}}^{\mathbb H}({\bf x})}\right|\,d{\bf x}\\
&=\,\dfrac{1}{\left|aa^{\prime}C_{\psi}\right|}\int_{\mathbb R^2} \left|\psi\left(\frac{{\bf x}-{\bf y}}{a}\right)\right|_{\mathbb H}\,\,\left|\overline{\psi\left(\frac{{\bf x}-{\bf y}^{\prime}}{a^{\prime}}\right)}\right|_{\mathbb H}\,d{\bf x}\\
&=\,\dfrac{1}{\left|aa^{\prime}C_{\psi}\right|}\int_{\mathbb R^2} \left|\psi\left(z-\frac{{\bf y}}{a}\right)\right|_{\mathbb H}\,\left|\psi\left(z-\frac{{\bf y}^{\prime}}{a^{\prime}}\right)\right|_{\mathbb H}\,aa^{\prime}dz\\
&=\,\dfrac{1}{\left|C_{\psi}\right|}\int_{\mathbb R^2} \left|\psi\left(z\right)\right|^2_{\mathbb H}\,dz\\
&=C_{\psi}^{-1}\,\big\|\psi\big\|_{L^2\left(\mathbb R^{2},\,\mathbb H\right)},\, \text{if}\, C_{\psi}>0.
\end{align*}

This completes the proof of theorem. \quad\fbox

\section{Uncertainty Principles for Quaternion Linear Canonical Wavelet Transform}

\parindent=0mm \vspace{.1in}
Heisenberg's uncertainty principle  in harmonic analysis is of central importance in time-frequency analysis as is it provides a lower bound for optimal simultaneous resolution in the time and frequency domains (see \cite{Fol}). This principle has been extended to different time-frequency transforms and several other versions of the uncertainty principle have been investigated from time to time. For instance, Beckner \cite{Bec} obtained a logarithmic version of the uncertainty principle by using a sharp form of Pitt's inequality and showed that this version yields the classical Heisenberg's inequality by virtue of Jensen's inequality. In this Section, we shall establish an analogue of the well-known Heisenberg's uncertainty inequality and the corresponding logarithmic version for the quaternion linear canonical wavelet transform as defined by (3.9). We first prove the following lemma.

\parindent=0mm \vspace{.1in}

{\bf Lemma 4.1.} {\it Let $\psi\in L^2(\mathbb R^2,\,\mathbb H)$ be an admissible quaternion wavelet, then for every $f \in L^2(\mathbb R^2,\,\mathbb H),$ we have }

\begin{align*}
 \int_{\mathcal G} \,\Big|w_k\mathcal L^{\mathbb H}_{A_1,A_2}\Big[\mathcal W_{\psi}^{\mathbb H}\big[ f\big]\Big]({\bf w})\Big|_{\mathbb H}^2 d\eta\,=\, C_{\psi}\left\|w_k \mathcal F_q[f]({\bf w})\right\|^2_{L^2(\mathbb R^2,\,\mathbb H)}\, \tag{4.1}
\end{align*}

\parindent=0mm \vspace{.0in}
{\it Proof.} Invoking Parseval formula (3.13) of quaternion linear canonical wavelet transform, we have

\begin{align*}
C_{\psi} \langle f,\,g\rangle_{L^2(\mathbb R^2,\mathbb H)} = \int_{\mathbb R^{+}}\int_{\text{SO(2)} }\int_{\mathbb R^2} \mathcal W_{\psi}^{\mathbb H}\big[ f\big](a,{\bf y},\theta)\,\overline{\mathcal W_{\psi}^{\mathbb H}\big[ g\big](a,{\bf y},\theta)}\,\dfrac{da d{\bf y} d\theta}{a^3}\\
\end{align*}

Implementing Plancheral theorem of QFT(two-sided) on L.H.S and Plancheral theorem of QLCT(two-sided) to the ${\bf y}$-integral on R.H.S of above equation, we obtain

\begin{align*}
C_{\psi} \left\langle \mathcal F_q[f],\,\mathcal F_q[g]\right\rangle_{L^2(\mathbb R^2,\mathbb H)} = \int_{\mathbb R^{+}}\int_{\text{SO(2)} }\left\langle\mathcal L^{\mathbb H}_{A_1,A_2}\Big[\mathcal W_{\psi}^{\mathbb H}\big[ f\big]\Big],\,\mathcal L^{\mathbb H}_{A_1,A_2}\Big[\mathcal W_{\psi}^{\mathbb H}\big[ g\big]\Big]\right\rangle_{L^2(\mathbb R^2,\mathbb H)}\,\dfrac{da d\theta}{a^3}\\
\end{align*}

Multiplying $\left|w_k\right|$ on both sides both sides of above equation, we get

\begin{align*}
C_{\psi} \left\langle w_k\,\mathcal F_q[f],\,w_k\,\mathcal F_q[g]\right\rangle_{L^2(\mathbb R^2,\mathbb H)} = \int_{\mathbb R^{+}}\int_{\text{SO(2)} }\left\langle w_k\,\mathcal L^{\mathbb H}_{A_1,A_2}\Big[\mathcal W_{\psi}^{\mathbb H}\big[ f\big]\Big],\,w_k\,\mathcal L^{\mathbb H}_{A_1,A_2}\Big[\mathcal W_{\psi}^{\mathbb H}\big[ g\big]\Big]\right\rangle_{L^2(\mathbb R^2,\mathbb H)}\,\dfrac{da d\theta}{a^3}\\
\end{align*}

Finally, for $f=g$, above equation yields

\begin{align*}
C_{\psi}\int_{\mathbb R^2} \left|w_k\,\mathcal F_q[f]({\bf w})\right|_{\mathbb H}^2 d{\bf w} =\int_{\mathbb R^{+}}\int_{\text{SO(2)} }\int_{\mathbb R^2} \,\Big|w_k\mathcal L^{\mathbb H}_{A_1,A_2}\Big[\mathcal W_{\psi}^{\mathbb H}\big[ f\big]\Big]({\bf w})\Big|_{\mathbb H}^2 \dfrac{da d{\bf w} d\theta}{a^3}\,
\end{align*}

This completes the proof of Lemma 4.1.\quad \fbox

\parindent=8mm \vspace{.1in}
We are now ready to establish the Heisenberg-type inequalities for the proposed quaternion linear canonical wavelet transform as defined by (3.9).

\parindent=0mm \vspace{.1in}
{\bf Theorem 4.2.} {\it Let $\psi\in L^2(\mathbb R^2,\,\mathbb H)$ be an admissible quaternion wavelet, the quaternion linear canonical wavelet transform $\mathcal W_{\psi}^{\mathbb H}\big[ f\big](a,{\bf y},\theta)$ given by (3.9) satisfies the following uncertainty inequality:}
\begin{align*}
\int_{\mathbb G} \left|y_k^2 \mathcal W_{\psi}^{\mathbb H}\big[ f\big](a,{\bf y},\theta)\right|_{\mathbb H}^2 d{\bf y} \cdot \int_{\mathbb R^2} \Big|w_k\,\mathcal F_q\big[f\big]({\bf w})\Big|_{\mathbb H}^2 d{\bf w} \ge \dfrac{b_k^2}{4}\sqrt{C_{\psi}}\big\|f\big\|^4_{L^2(\mathbb R^2,\,\mathbb H)} \tag{4.2}
\end{align*}

\parindent=0mm \vspace{.1in}
{\it Proof.} Invoking the Heisenberg's  inequality for the quaternion linear canonical transform \cite{Kou}, we can write
\begin{align*}
\left\{\int_{\mathbb R^2} {y_k}^2\big|f({\bf y})\big|_{\mathbb H}^2d{\bf y} \right\}^{1/2} \left\{\int_{\mathbb R^2} w_k^2 \Big|\mathcal{L}_{A_1,A_2}^{\mathbb H} \left[f\right]({\bf w})\Big|_{\mathbb H}^2 d{\bf w}\right\}^{1/2}\ge \dfrac{b_k}{2} \int_{\mathbb R^2}\big|f({\bf x})\big|_{\mathbb H}^2 d{\bf x}.
\end{align*}

Integrating both sides of the above inequality with respect to measure $\frac{da\,d\theta}{a^3}$,\, we have
\begin{align*}
&\int_{\mathbb R^{+}}\int_{\text{SO(2)}}\Bigg[ \left\{\int_{\mathbb R^2} {y_k}^2\big|f({\bf y})\big|_{\mathbb H}^2d{\bf y} \right\}^{1/2} \left\{\int_{\mathbb R^2} w_k^2 \Big|\mathcal{L}_{A_1,A_2}^{\mathbb H} \left[f\right]({\bf w})\Big|_{\mathbb H}^2 d{\bf w}\right\}^{1/2}\Bigg] \dfrac{da\,d\theta}{a^3}\qquad\qquad\\
&\qquad\qquad\qquad\qquad\qquad\qquad\qquad\qquad\qquad \ge \dfrac{b_k}{2} \int_{\mathbb R^{+}}\int_{\text{SO(2)}} \int_{\mathbb R^2}\big|f({\bf x})\big|_{\mathbb H}^2 \dfrac{da\,d{\bf x}\,d\theta}{a^3}
\end{align*}

Now by implementing Cauchy-Schwartz inequality, we obtain
\begin{align*}
\left\{\int_{\mathbb R^{+}}\int_{\text{SO(2)}}\int_{\mathbb R^2} {y_k}^2\big|f({\bf y})\big|_{\mathbb H}^2d{\bf y}\dfrac{da\,d\theta}{a^3} \right\}^{1/2} &\left\{\int_{\mathbb R^{+}}\int_{\text{SO(2)}}\int_{\mathbb R^2} w_k^2 \Big|\mathcal{L}_{A_1,A_2}^{\mathbb H} \left[f\right]({\bf w})\Big|_{\mathbb H}^2 d{\bf w}\dfrac{da\,d\theta}{a^3}\right\}^{1/2}\\
&\qquad\qquad\ge \dfrac{b_k}{2} \int_{\mathbb R^{+}}\int_{\text{SO(2)}} \int_{\mathbb R^2}\big|f({\bf x})\big|_{\mathbb H}^2 \dfrac{da\,d{\bf x}\,d\theta}{a^3}
\end{align*}

Considering   $\mathcal W_{\psi}^{\mathbb H}\big[ f\big](a,{\bf y},\theta)$  as a function of ${\bf y}$  and replacing $f$ by  $\mathcal W_{\psi}^{\mathbb H}\big[ f\big](a,{\bf y},\theta)$ in above equation to obtain

\begin{align*}
&\left\{\int_{\mathbb R^{+}}\int_{\text{SO(2)}}\int_{\mathbb R^2} {y_k}^2\left|\mathcal W_{\psi}^{\mathbb H}\big[ f\big](a,{\bf y},\theta)\right|_{\mathbb H}^2d{\bf y}\dfrac{da\,d\theta}{a^3} \right\}^{1/2}\\
&\qquad\times\left\{\int_{\mathbb R^{+}}\int_{\text{SO(2)}}\int_{\mathbb R^2} w_k^2 \Big|\mathcal{L}_{A_1,A_2}^{\mathbb H} \left[\mathcal W_{\psi}^{\mathbb H}\big[ f\big](a,{\bf y},\theta)\right]({\bf w})\Big|_{\mathbb H}^2 d{\bf w}\dfrac{da\,d\theta}{a^3}\right\}^{1/2}\\
&\qquad\qquad\qquad\qquad\qquad\qquad\qquad\ge \dfrac{b_k}{2} \int_{\mathbb R^{+}}\int_{\text{SO(2)}} \int_{\mathbb R^2}\big|\mathcal W_{\psi}^{\mathbb H}\big[ f\big](a,{\bf y},\theta)\big|_{\mathbb H}^2 \dfrac{da\,d{\bf y}\,d\theta}{a^3}
\end{align*}

By virtue of Lemma 4.1,  we obtain

\begin{align*}
&\left\{\int_{\mathbb R^{+}}\int_{\text{SO(2)}}\int_{\mathbb R^2} {y_k}^2\left|\mathcal W_{\psi}^{\mathbb H}\big[ f\big](a,{\bf y},\theta)\right|_{\mathbb H}^2d{\bf y}\dfrac{da\,d\theta}{a^3} \right\}^{1/2}\left\{C_{\psi}\int_{\mathbb R^2} \left|w_k\,\mathcal F_q[f]({\bf w})\right|_{\mathbb H}^2 d{\bf w}\right\}^{1/2}\\
&\qquad\qquad\qquad\qquad\qquad\qquad\qquad\qquad\qquad\ge \dfrac{b_k}{2} \int_{\mathbb R^{+}}\int_{\text{SO(2)}} \int_{\mathbb R^2}\big|\mathcal W_{\psi}^{\mathbb H}\big[ f\big](a,{\bf y},\theta)\big|_{\mathbb H}^2 \dfrac{da\,d{\bf y}\,d\theta}{a^3}
\end{align*}

Furthermore, on implementing Corollary 3.10 on R.H.S of the above inequality we have

\begin{align*}
&\left\{\int_{\mathbb R^{+}}\int_{\text{SO(2)}}\int_{\mathbb R^2} {y_k}^2\left|\mathcal W_{\psi}^{\mathbb H}\big[ f\big](a,{\bf y},\theta)\right|_{\mathbb H}^2d{\bf y}\dfrac{da\,d\theta}{a^3} \right\}^{1/2}\left\{C_{\psi}\int_{\mathbb R^2} \left|w_k\,\mathcal F_q[f]({\bf w})\right|_{\mathbb H}^2 d{\bf w}\right\}^{1/2}\\
&\qquad\qquad\qquad\qquad\qquad\qquad\qquad\qquad\qquad\ge \dfrac{b_k}{2} C_{\psi}\, \big\|f\big\|^{2}_{L^2(\mathbb R^2,\mathbb H)}
\end{align*}

Dividing both sides by $\sqrt{C_{\psi}}$, we obtain the desired result.\quad \fbox

\parindent=8mm \vspace{.1in}
Before presenting our next result, we have the following definition of space of rapidly decreasing smooth quaternion functions(see \cite{Asym}).

\parindent=0mm \vspace{.1in}
{\bf Definition 4.3.} For a multi-index $\alpha =(\alpha_1,\alpha_2)\in \mathbb R^{+}\times \mathbb R^{+},$  the Schwartz space in $ L^2(\mathbb R^2,\mathbb H)$ is defined as

\begin{align*}
\mathcal S (\mathbb R^2,\mathbb H)= \left\{f\in\mathbb C^{\infty} (\mathbb R^2,\mathbb H); \sup_{t\in\mathbb R^2}\big(1+|t|^k\big)\Big|\dfrac{\partial^{\alpha_1+\alpha_2}[f(t)]}{\partial^{\alpha_1}_{t_1}\partial^{\alpha_2}_{t_2}}\Big|\,< \infty \right\},
\end{align*}
where $\mathbb C^{\infty}(\mathbb R^2,\mathbb H)$ is the set of smooth functions from $\mathbb R^2$ to $\mathbb H$.

\parindent=8mm \vspace{.1in}

We now establish the logarithmic uncertainty principle for the quaternion linear canonical wavelet transform $\mathcal W_{\psi}^{\mathbb H}\big[ f\big]$ as defined by (3.9).

\parindent=0mm \vspace{.1in}

{\bf Theorem 4.4.}{\it For an admissible quaternion wavelet $\psi\in \mathcal S (\mathbb R^2,\mathbb H)$ and a signal $f\in \mathcal S(\mathbb R^2,\mathbb H)$, the quaternion linear canonical wavelet transform $\mathcal W_{\psi}^{\mathbb H}\big[ f\big]$ satisfies the following logarithmic estimate of the uncertainty inequality:}
\begin{align*}
&\int_{\mathcal G}\ln|{\bf y}|\, \Big|\mathcal W_{\psi}^{\mathbb H}\big[ g\big](a,{\bf y},\theta)\Big|_{\mathbb H}^{2} d\eta +C_{\psi}\int_{\mathbb R^2}\ln|{\bf w}|\,\left|\mathcal F_q\big[f\big]({\bf w})\right|_{\mathbb H}^{2}d{\bf w}\,\geq\,C_{\psi}\big(\mathcal D+ln|b|\big)\int_{\mathbb R^2}\Big|f({\bf x})\Big|^2_{\mathbb H}d{\bf x},\, \tag{4.3}
\end{align*}
{\it where $\mathcal D=\left(\dfrac{\Gamma^{\prime}(1/2)}{\Gamma(1/2)}-\ln \pi\right)$ and $\Gamma$ is a Gamma function.}

\parindent=0mm \vspace{.1in}
{\it Proof.} For the quaternion-valued function $f\in L^2(\mathbb R^2,\mathbb H)$, the time and frequency spreads satisfy the inequality \cite{Bahri}
\begin{align*}
\int_{\mathbb R^{2}}\ln|{\bf y}|\,\Big|f({\bf y})\Big|_{\mathbb H}^{2}d{\bf y}+\int_{\mathbb R^{2}}\ln|{\bf w}|\,{\Big|\mathcal{L}_{A_1,A_2}^{\mathbb H} \big[f\big]({\bf w})\Big|_{\mathbb H}^{2}}\,d{\bf w}\geq\,\big(\mathcal D+ln|b|\big)\int_{\mathbb R^{2}}\Big|f({\bf y})\Big|_{\mathbb H}^{2}\,d{\bf y}
\end{align*}

Replacing $f({\bf y})$ by  $\mathcal W_{\psi}^{\mathbb H}\big[ f\big](a,{\bf y},\theta)$ in the above inequality, we obtain

\begin{align*}
&\int_{\mathbb R^{2}}\ln|{\bf y}|\,\left|\mathcal W_{\psi}^{\mathbb H}\big[ f\big](a,{\bf y},\theta)\right|^2_{\mathbb H} d{\bf y}+\int_{\mathbb R^{2}}\ln|{\bf w}|\,\left|\mathcal{L}_{A_1,A_2}^{\mathbb H} \left[\mathcal W_{\psi}^{\mathbb H}\big[ f\big]\right]({\bf w})\right|_{\mathbb H}^{2}\,d{\bf w}\\
&\qquad\qquad\qquad\qquad\qquad\qquad\qquad\qquad\qquad \geq\,\big(\mathcal D+ln|b|\big)\int_{\mathbb R^{2}}\left|\mathcal W_{\psi}^{\mathbb H}\big[ f\big](a,{\bf y},\theta)\right|_{\mathbb H}^{2}\,d{\bf y}
\end{align*}

Integrating above equation with respect to measure $\frac{dad\theta}{a^3}$, and then applying the  Fubini theorem, we obtain
\begin{align*}
&\int_{\mathbb R^{+}}\int_{\text{SO(2)}}\int_{\mathbb R^{2}}\ln|{\bf y}|\,\left|\mathcal W_{\psi}^{\mathbb H}\big[ f\big](a,{\bf y},\theta)\right|^2_{\mathbb H} d{\bf y}\dfrac{dad\theta}{a^3}+\int_{\mathbb R^{+}}\int_{\text{SO(2)}}\int_{\mathbb R^{2}}\ln|{\bf w}|\,\left|\mathcal{L}_{A_1,A_2}^{\mathbb H} \left[\mathcal W_{\psi}^{\mathbb H}\big[ f\big]\right]({\bf w})\right|_{\mathbb H}^{2}\,d{\bf w}\,\dfrac{dad\theta}{a^3}\\
&\qquad\qquad\qquad\qquad\qquad\qquad\qquad\qquad\geq\,\big(\mathcal D+ln|b|\big)\int_{\mathbb R^{+}}\int_{\text{SO(2)}}\int_{\mathbb R^{2}}\left|\mathcal W_{\psi}^{\mathbb H}\big[ f\big](a,{\bf y},\theta)\right|_{\mathbb H}^{2}\,d{\bf y}\,\dfrac{dad\theta}{a^3}
\end{align*}

Applying Lemma 4.1 for $w_k^2=\ln|w_k|$ on L.H.S and Corollary 3.10 on R.H.S, we obtain the desired result
\begin{align*}
&\int_{\mathcal G}\ln|{\bf y}|\, \Big|\mathcal W_{\psi}^{\mathbb H}\big[ g\big](a,{\bf y},\theta)\Big|_{\mathbb H}^{2} d\eta +C_{\psi}\int_{\mathbb R^2}\ln|{\bf w}|\,\left|\mathcal F_q\big[f\big]({\bf w})\right|_{\mathbb H}^{2}d{\bf w}\,\geq\,C_{\psi}\big(\mathcal D+ln|b|\big)\int_{\mathbb R^2}\Big|f({\bf x})\Big|^2_{\mathbb H}d{\bf x},\,
\end{align*}

This completes the proof of Theorem 4.4. \quad\fbox

In the following, we establish a local type uncertainty principle for quaternion linear canonical wavelet transform $\mathcal W_{\psi}^{\mathbb H}\big[ f\big]$ as defined by (3.9).

\parindent=0mm \vspace{.1in}

{\bf Theorem 4.5.} {\it Given an admissible quaternion wavelet $\psi\in L^2 (\mathbb R^2,\mathbb H)$ and a signal $f\in L^2(\mathbb R^2,\mathbb H)$, with $\|\psi\|^2_{L^2(\mathbb R^2,\mathbb H)}=1=\|f\|^2_{L^2(\mathbb R^2,\mathbb H)},$
such that for a measurable set $E\subset \mathbb R^2\times\mathbb R^2,\, \epsilon\geq 0,$ and}
\begin{align*}
\int \int_{E} \Big|\mathcal W_{\psi}^{\mathbb H}\big[ f\big](a,{\bf y},\theta)\Big|_{\mathbb H}^{2} d{\bf y} d{\bf x}\,\geq\,1-\epsilon.\tag{4.4}
\end{align*}
We have $a (1-\epsilon)\leq \,\mu(E),$ where $\mu(E)$ is Lebesgue measure of $E$.

\parindent=0mm \vspace{.1in}
{\it Proof.} From Definition 3.4, we have
\begin{align*}
\Big|\mathcal W_{\psi}^{\mathbb H}\big[ f\big](a,{\bf y},\theta) \Big|_{\mathbb H}&=\Big|  a^{-1} \int_{\mathbb R^{2}} f({\bf x})\,e^{\frac{ja_2}{2b_2}(x_2^2-y_2^2)}\, \overline{\psi\big(r_{-\theta}a^{-1}({\bf x}-{\bf y})\big)}\,e^{\frac{ia_1}{2b_1}(x_1^2-y_1^2)} \,d{\bf x}\Big|_{\mathbb H}
\\&\leq\,\dfrac{1}{a}\int_{\mathbb R^2}|f({\bf x})|\,\big|\overline{\psi\big(r_{-\theta}a^{-1}({\bf x}-{\bf y})\big)}\big|\,d{\bf x}.
\end{align*}

Taking Sup-norm on L.H.S and implementing well known Holders inequality on R.H.S, we obtain

\begin{align*}
\Big\|\mathcal W_{\psi}^{\mathbb H}\big[ f\big](a,{\bf y},\theta)\Big\|_{L^{\infty}(\mathbb R^2,\,\mathbb H)}
\leq\,\dfrac{1}{a}\,\big\|f\big\|_{L^{2}(\mathbb R^2,\,\mathbb H)}\big\|\psi\big\|_{L^{2}(\mathbb R^2,\,\mathbb H)}\,\tag{4.5}
\end{align*}

Plugging inequality (4.5) in (4.4), we get
\begin{align*}
1-\epsilon\,&\leq\,\int \int_{E} \Big|\mathcal W_{\psi}^{\mathbb H}\big[ f\big](a,{\bf y},\theta)\Big|_{\mathbb H}^{2} d{\bf y} d{\bf x}\\
&\leq\, \mu(E)\cdot \Big\|\mathcal W_{\psi}^{\mathbb H}\big[ f\big](a,{\bf y},\theta) \Big\|_{L^{\infty}(\mathbb R^2,\,\mathbb H)}\\
&\leq\, \dfrac{\mu(E)}{a}\,\big\|f\big\|_{L^{2}(\mathbb R^2,\,\mathbb H)}\big\|\psi\big\|_{L^{2}(\mathbb R^2,\,\mathbb H)}\\
&\leq\, \dfrac{\mu(E)}{a}.
\end{align*}

This completes the proof of Theorem 4.5. \quad\fbox

\parindent=0mm \vspace{.1in}

{\bf Theorem 4.6(Local uncertainty inequality).} {\it Let $E$ be a measurable subset of $\mathbb R^2\times\mathbb R^2$, such that $0<\mu(E)<1,$ then for every
$f\in L(\mathbb R^2,\mathbb H)$ and an admissible quaternion wavelet $\psi\in L (\mathbb R^2,\mathbb H)$ , we have}

\begin{align*}
\big\|f\big\|_{L^{2}(\mathbb R,\,\mathbb H)}\big\|\psi\big\|_{L^{2}(\mathbb R,\,\mathbb H)}\leq \dfrac{1}{\sqrt{1-\mu(E)}}
\Big\|\mathcal W_{\psi}^{\mathbb H}\big[ f\big](a,{\bf y},\theta)\Big\|_{L^2(E^c,\,\mathbb H)}\,\tag{4.6}
\end{align*}

{\it Moreover, for every $\alpha>0,$ there exist $C(\alpha)>0,$ such that}

\begin{align*}
\big\|f\big\|_{L^{2}(\mathbb R,\,\mathbb H)}\big\|\psi\big\|_{L^{2}(\mathbb R,\,\mathbb H)}\leq C(\alpha)\,\Big[\int_{\mathbb R^2}\int_{\mathbb R^2}\big|({\bf y},{\bf x})\big|^{2\alpha}\,
\Big|\mathcal W_{\psi}^{\mathbb H}\big[ f\big](a,{\bf y},\theta)\Big|^2_{\mathbb H}\,d{\bf y} d{\bf x}\Big]^{1/2}.\,\tag{4.7}
\end{align*}

\parindent=0mm \vspace{.1in}
{\it Proof.} By invoking the Corollary 3.10 and Theorem 4.5, we have

\begin{align*}
\Big\|\mathcal W_{\psi}^{\mathbb H}\big[ f\big](a,{\bf y},\theta)\Big\|^2_{L^2(\mathbb R^2\times\mathbb R^2,\,\mathbb H)}
&=\int_{\mathbb R^2}\int_{\mathbb R^2} \Big|\mathcal W_{\psi}^{\mathbb H}\big[ f\big](a,{\bf y},\theta)\Big|^2_{\mathbb H}d{\bf y} d{\bf x}\\
&=\int \int_{E} \Big|\mathcal W_{\psi}^{\mathbb H}\big[ f\big](a,{\bf y},\theta)\Big|^2_{\mathbb H}d{\bf y} d{\bf x}+\int \int_{E^c} \Big|\mathcal W_{\psi}^{\mathbb H}\big[ f\big](a,{\bf y},\theta)\Big|^2_{\mathbb H}d{\bf y} d{\bf x}\\
&\leq \,\mu(E) \big\|f\big\|^2_{L^{2}(\mathbb R,\,\mathbb H)}\big\|\psi\big\|^2_{L^{2}(\mathbb R,\,\mathbb H)}+\Big\|\mathcal W_{\psi}^{\mathbb H}\big[ f\big](a,{\bf y},\theta)\Big\|^2_{L^2(E^c,\,\mathbb H)}
\end{align*}

Equivalently,
\begin{align*}
(1-\mu(E))\big\|f\big\|^2_{L^{2}(\mathbb R,\,\mathbb H)}\big\|\psi\big\|^2_{L^{2}(\mathbb R,\,\mathbb H)}\leq \Big\|\mathcal W_{\psi}^{\mathbb H}\big[ f\big](a,{\bf y},\theta)\Big\|^2_{L^2(E^c,\,\mathbb H)}
\end{align*}

Taking square root on both sides and then dividing both sides by $\sqrt{1-\mu(E)},$ we get

\begin{align*}
\big\|f\big\|_{L^{2}(\mathbb R,\,\mathbb H)}\big\|\psi\big\|_{L^{2}(\mathbb R,\,\mathbb H)}\leq \dfrac{1}{\sqrt{1-\mu(E)}}\Big\|\mathcal W_{\psi}^{\mathbb H}\big[ f\big](a,{\bf y},\theta)\Big\|_{L^2(E^c,\,\mathbb H)}
\end{align*}

which proves our first assertion.

\parindent=0mm \vspace{.0in}
Now, we fix ${\lambda}_{0}\in(0,1]$ small enough such that $\mu\big(B_{{\lambda}_{0}}\big)<1,$ where $B_{{\lambda}_{0}}=\left\{({\bf y},{\bf x})\in \mathbb R^2\times\mathbb R^2;\,\big|({\bf y},{\bf x})\big|<\,{\lambda}_{0} \right\},$
the ball of radius ${\lambda}_{0}$ centered at origin, we have from inequality (4.6),

\begin{align*}
\big\|f\big\|_{L^{2}(\mathbb R,\,\mathbb H)}\big\|\psi\big\|_{L^{2}(\mathbb R,\,\mathbb H)}
&\leq\, \dfrac{1}{\sqrt{1-\mu(B_{{\lambda}_{0}})}}\left[\int \int_{\big|({\bf y},{\bf x})\big|>\,{\lambda}_{0} }\Big|\mathcal W_{\psi}^{\mathbb H}\big[ f\big](a,{\bf y},\theta)\Big|^2_{\mathbb H}d{\bf y} d{\bf x} \right]^{1/2}\\
&\leq\, \dfrac{1}{{\lambda}^{\alpha}_0 \sqrt{1-\mu(B_{{\lambda}_{0}})}}\left[\int \int_{\big|({\bf y},{\bf x})\big|>\,{\lambda}_{0} }{\lambda}^{2\alpha}_0\Big|\mathcal W_{\psi}^{\mathbb H}\big[ f\big](a,{\bf y},\theta)\Big|^2_{\mathbb H}d{\bf y} d{\bf x} \right]^{1/2}\\
&\leq\, \dfrac{1}{{\lambda}^{\alpha}_0 \sqrt{1-\mu(B_{{\lambda}_{0}})}}\left[\int \int_{\big|({\bf y},{\bf x})\big|>\,{\lambda}_{0}} \,\big|({\bf y},{\bf x})\big|^{2\alpha}\Big|\mathcal W_{\psi}^{\mathbb H}\big[ f\big](a,{\bf y},\theta)\Big|^2_{\mathbb H}d{\bf y} d{\bf x} \right]^{1/2}\\
&\leq\, C(\alpha)\left[\int \int_{\big|({\bf y},{\bf x})\big|>\,{\lambda}_{0}} \,\big|({\bf y},{\bf x})\big|^{2\alpha}\Big|\mathcal W_{\psi}^{\mathbb H}\big[ f\big](a,{\bf y},\theta)\Big|^2_{\mathbb H}d{\bf y} d{\bf x} \right]^{1/2}
\end{align*}
where $C(\alpha)=\dfrac{1}{{\lambda}^{\alpha}_{0} \sqrt{1-\mu(B_{{\lambda}_{0}})}}.$ This completes the proof of  Theorem 4.6. \quad\fbox

\end{document}